\documentclass[11pt,notitlepage]{amsart}

\newcommand{\E}{\mathbb{E}}
\newcommand{\Prob}{\mathbb{P}}
\newcommand{\R}{\mathbb{R}}

\newcommand{\ud}{\mathrm{d}}
\newcommand{\udx}{\mathrm{dx}}
\newcommand{\udy}{\mathrm{dy}}
\newcommand{\udz}{\mathrm{dz}}
\newcommand{\ai}{\mathrm{Ai}}
\newcommand{\dai}{\mathrm{Ai'}}
\newcommand{\var}{\mathrm{Var}}

\newtheorem{theorem}{Theorem}[section]
\newtheorem{lemma}{Lemma}[section]

\newenvironment{axiom-T}{\hspace{-15pt}\textbf{Axiom-T} } {$\;$}  
\newenvironment{axiom-S}{\hspace{-15pt}\textbf{Axiom-S} } {$\;$} 
\newenvironment{axiom-M}{\hspace{-15pt}\textbf{Axiom-M} } {$\;$}  
\newenvironment{axiom-PH}{\hspace{-15pt}\textbf{Axiom-PH} } {$\;$} 
\newenvironment{axiom-R}{\hspace{-15pt}\textbf{Axiom-R} } {$\;$}

\usepackage{epsfig}
\usepackage{amsmath,amssymb}

\title{Gaussian fluctuations of eigenvalues in the GUE}
\author[J. Gustavsson]{Jonas Gustavsson}
\address{Jonas Gustavsson \newline \indent Royal institute of Technology}
\email{jonasg@math.kth.se}

\begin{document}

\begin{abstract}
Under certain conditions on $k$ we calculate the limit distribution of 
eigenvalue number $k$, $x_k$, of the Gaussian unitary 
ensemble (GUE). More specifically, if $n$ is the dimension of a random matrix
from the GUE and $k$ is such that both $k$ and $n-k$ tend to infinity as 
$n \rightarrow \infty$, then $x_k$ is normally distributed in the limit.

We also consider the joint limit distribution of $(x_{k_1}, \ldots ,x_{k_m})$
where we require that $k_1$, $n-k_m$ and $k_{i+1}-k_i$, $1 \leq i \leq m-1$, 
tend to infinity with $n$. 
The result is an $m$-dimensional normal distribution.
\end{abstract}

\maketitle

\section{Introduction and formulation of results}
\setcounter{equation}{0}

The Gaussian unitary ensemble (GUE) is a classical random matrix ensemble. 
It is defined by the probability distribution on the space of n$\times$n
Hermitian matrices given by
\begin{equation*}
\Prob(dH) = C_n \, e^{-\textrm{Trace}H^2} dH.
\end{equation*}
By $dH$ we mean the Lebesgue measure on the $n^2$ essentially different
members of the matrix, namely
\begin{equation} \label{members}
\{ \textrm{Re} H_{ij}; 1 \leq i \leq j \leq n, 
\textrm{Im} H_{ij}; 1 \leq i < j \leq n \} .
\end{equation}
In other words this means that the entries in (\ref{members}) are
independent Gaussian random variables with zero mean and variance 
$\frac{1+\delta_{ij}}{4}$.
The measure on the matrices naturally induces a measure on the 
corresponding $n$ real eigenvalues $x_i$. This induced measure can be
explicitly calculated and its density is given by
\begin{equation*}
p_n(x_1, \ldots, x_n) =
\frac{1}{Z^{(2)}_n} \prod_{1 \leq i < j \leq n}
|x_i-x_j|^2 \, e^{-x_1^2 - \ldots -x_n^2} .
\end{equation*}
The normalization constant $Z_n^{(2)}$ is called the partition function. 
It is often convenient to work with the eigenvalues being ordered.
Naming the eigenvalues so that $x_1 < \ldots < x_n$, gives that the 
probability density
$\rho_{n,n}(x_1, \ldots ,x_n)$ of the ordered eigenvalues defined on the space 
\begin{equation*}
\R^n_{\mathrm{ord}} = \left\{ x_1, \ldots ,x_n ; x_1 < \ldots <x_n \right\}
\end{equation*}
is given by 
\begin{equation*}
\rho_{n,n}(x_1, \ldots ,x_n) = n! \, p_n(x_1, \ldots, x_n) .
\end{equation*}
This density $\rho$ is a member of a family of functions called the
correlation functions. These functions are defined by
\begin{multline*}
\rho_{n,k}(x_1, \ldots,x_k) =
\frac{n!}{(n-k)!} \int_{\mathbb{R}^{n-k}} p_n(x_1, \ldots, x_n) \,
\udx_{k+1} \ldots \udx_n \\
= \mathrm{det} (K_n(x_i,x_j))_{i,j=1}^k .
\end{multline*}
Here $K_n(x,y)$ is given by
\begin{equation*}
K_n(x,y) = \sum_{i=0}^{n-1} h_i(x) h_i(y) e^{- \frac{1}{2}(x^2 + y^2)}
\end{equation*}
where $\{ h_i \}_{i \geq 0}$ are the orthonormalized 
Hermite polynomials, that is
\begin{equation*}
\int_{-\infty}^{\infty} h_i(x) h_j(x) e^{-x^2} \, \udx = \delta_{ij}.
\end{equation*}
The kernel
$K_n(x,y)$ can also be represented by the so called Christoffel-Darboux
identity. For $x \neq y$ it holds that
\begin{equation*}
K_n(x,y) = \left( \frac{n}{2} \right)^{1/2} 
\frac{h_n(x) h_{n-1}(y) - h_n(y) h_{n-1}(x)}{x-y}
e^{-\frac{1}{2}(x^2 + y^2)}
\end{equation*} 
and on the diagonal one has
\begin{equation*}
K_n(x,x) = \left( n h^2_n(x) - \sqrt{n(n+1)} h_{n-1}(x) h_{n+1}(x) \right)
e^{- x^2} .
\end{equation*}
The correlation function $\rho_{n,1}$ describes the overall density of 
the eigenvalues and the Wigner semi-circle law states that
\begin{equation} \label{semi}
\lim_{n \rightarrow \infty} \sqrt{\frac{2}{n}} \rho_{n,1} 
( \sqrt{2n}x ) =
\left\{ \begin{array}{ll}
\frac{2}{\pi} \sqrt{1-x^2} & \textrm{if $|x| \leq 1$} \\
0 & \textrm{if $|x| > 1$.} 
\end{array} \right.
\end{equation} 
All the results above and more can be found in the book by Mehta, \cite{me}. 

This paper deals with the distribution of eigenvalue number $k$, $x_k$, of the
GUE. More specifically, we look at the distribution of 
$x_k$ as $n$ and $k$ both tend to infinity. For example, if 
\begin{equation*}
k = k(n) = n - \log n
\end{equation*}
then, as $n$ becomes large, $k$ is very close (relatively) to the right
edge of the spectrum. Another example is when $k$ = $n/2$. In this case 
we are in the middle of the bulk of the spectrum. In both cases one 
ends up with a normal distribution in the limit. 
The following theorems generalize and specify this statement. 
\begin{theorem}[The bulk] \label{th1}
Set
\begin{equation*}
G(t) = \frac{2}{\pi} \int_{-1}^{t} \sqrt{1-x^2} \, \udx \qquad
-1 \leq t \leq 1 
\end{equation*}
and $t = t(k,n) = G^{-1}(k/n)$
where $k = k(n)$ is such that $k/n \rightarrow a \in (0,1)$
as $n \rightarrow \infty$. If $x_{k}$ denotes eigenvalue number $k$
in the GUE it holds that, as $n \rightarrow \infty$,
\begin{equation*}
\frac{x_k - t \sqrt{2n}}{\left( \frac{\log n}{4(1-t^2)n}
\right)^{1/2}} \longrightarrow \mathrm{N(0,1)}
\end{equation*}
in distribution.
\end{theorem}
 
\begin{theorem} [The edge] \label{th2}
Let $k$ be such that $k \rightarrow \infty$ but 
$\frac{k}{n} \rightarrow 0$ as $n \rightarrow \infty$ and 
let $x_{n-k}$ denote eigenvalue number $n-k$ in the GUE. 
Then it holds that, as $n \rightarrow \infty$, 
\begin{equation*}
\frac{x_{n-k} - \sqrt{2n} \left( 1 - 
\left( \frac{3 \pi k}{4 \sqrt{2} n} \right)^{2/3} \right)}{\left( 
\left( \frac{1}{12 \pi} \right)^{2/3} \frac{\log k}{n^{1/3} k^{2/3}}
\right)^{1/2}} \longrightarrow \mathrm{N(0,1)}
\end{equation*}
in distribution. 
\end{theorem}

\emph{Remark 1.} The theorems deal with the bulk and the right spectrum edge. 
One gets the equivalent for the left edge with some obvious modifications. \\  

\emph{Remark 2.} In \cite{tw} the distribution of the largest eigenvalue
was studied. \\

\emph{Remark 3.} Set $I = ( - \infty, s]$. In \cite{bs} it is shown that
\begin{multline} 
\Prob [x_k \in (s,s+\textrm{ds})] \\ = \label{Ja}
\left( \frac{1}{(k-1)!} \int_{I^{k-1}} 
J_k(x_1, \ldots, x_{k-1},s) \, \mu(\udx_1) \ldots \mu(\udx_{k-1}) \right)
\mu(\textrm{ds}).
\end{multline}
Here $J_k$ is the so called Janossy density and 
\begin{equation*}
\mu(\udx) = C \, e^{-\frac{x^2}{2}} \udx.
\end{equation*}
In \cite{bs} it is also proven that $J_k$ can be expressed explicitly by a 
determinantal formula.  
For $k$ and $n$ as in (\ref{th1}) or (\ref{th2}) we should thus have that
(\ref{Ja}) is, for large $n$, approximately equal to
the probability density function for the normal distribution $N(\nu,\sigma)$.
The parameters $\nu$ and $\sigma$ should of course be taken to be those 
indicated from the relevant theorem above. \\

\emph{Remark 4.} The zero number $k$ of the Hermite polynomial of degree
$n$ is close to the expected value of eigenvalue number $k$ of 
$\textrm{GUE}_n$. This can be shown directly by the following result, 
\cite{dkmvz}: 

There are constants $k_0$ and $C$ such that for 
$k_0 \leq k \leq n-k_0$ and $\alpha = k/n$ it holds that
\begin{multline}\label{zeroG}
\left| \frac{z_{k,n}}{\sqrt{2n}} -
G^{-1} \left[ \frac{k}{n} - \frac{1}{2 \pi n} 
\arcsin \left( G^{-1} (k/n) \right) + \frac{1}{2n} \right] \right| \\ \leq 
\frac{\textrm{C}}{n^2 (\alpha (1-\alpha))^{4/3}}.
\end{multline} 
Here $z_{1,n} < \ldots < z_{n,n}$ are the zeros of the Hermite polynomial 
of degree $n$. When we are in the bulk this translates into
\begin{equation*}
\left| z_{k,n} - \sqrt{2n} \, G^{-1}(k/n) \right| \leq \frac{C}{\sqrt{n}}.
\end{equation*}
This means that one can replace $t \sqrt{2n}$ by $z_{k,n}$ in Theorem
\ref{th1}. Close to the edge, for example when $k = n-\log{n}$, we can
not use (\ref{zeroG}) to conclude that this replacement is allowed. 

A motivation for this approximate equality between the locations of the
zeros and eigenvalues goes as follows. Set 
\begin{equation*}
W = \frac{1}{2} \sum_{i=1}^{n} x_i^2 - \sum_{1 \leq i<j\leq n}
\log |x_i - x_j| .
\end{equation*}
and note that 
\begin{equation*}
\rho_{n,n}(x_1, \ldots,x_n) = \textrm{Const} \cdot
e^{-2 W} .
\end{equation*}
It is a fact, see \cite{me}, that $W$ obtains its minimum exactly when 
$x_i = z_{i,n}$, $1 \leq i \leq n$. This configuration is hence the 
most 'probable' for the eigenvalues.
Expanding around this minimum we see that it is reasonable that
$x_k$ should have Gaussian fluctuations around $z_{k,n}$.\\

\emph{Remark 5.} If one is interested in the distribution of the 
eigenvalues of some other ensemble one should in many cases be able 
to apply the same methodology that we use in this paper. \\\\
It is also interesting to see what happens when looking at several
eigenvalues at the same time. When we write $k(n) \sim n^{\theta}$ 
below we mean
that $k(n) = h(n) n^{\theta}$ where $h$ is a function such that, 
for all $\epsilon > 0$,
\begin{equation} \label{sim}
\frac{h(n)}{n^{\epsilon}} \rightarrow 0 \quad \textrm{and} \quad
h(n) n^{\epsilon} \rightarrow \infty
\end{equation}
as $n \rightarrow \infty$. We have the following results:
\begin{theorem} [The bulk] \label{th3}
Let $\{x_{k_i}\}_{i=1}^{m}$ be eigenvalues of the GUE such that 
$0 < k_{i} - k_{i+1} \sim n^{\theta_i}$, $0 < \theta_i \leq1$, and 
$k_i/n \rightarrow a_i \in (0 ,1)$ 
as $n \rightarrow \infty$. Define $s_i=s_i(k_i,n)=G^{-1}(k_i/n)$ and
set
\begin{equation*}
X_i = \frac{x_{k_i} - s_i \sqrt{2n}}
{\left( \frac{\log n}{4(1-s_i^2)n}\right)^{1/2}}
\quad i = 1, \ldots , m.
\end{equation*}
Then as, $n \rightarrow \infty$, 
\begin{equation*}
\mathbb{P} \left[ X_1 \leq x_1, \ldots , X_m \leq x_m \right] 
\longrightarrow \Phi_{\Lambda}(x_1, \ldots ,x_m) 
\end{equation*}
where $\Phi_{\Lambda}$ is the cdf 
\footnote{Cumulative distribution function} for the 
$m$-dimensional normal distribution with covariance matrix 
$\Lambda_{i,j}= 1 - \max \{\theta_k | \, i \leq k < j < m \}$ if $i<j$ and
$\Lambda_{i,i} =1$.
\end{theorem}

\begin{theorem} [The edge] \label{th4}
Let $\{x_{n-k_i}\}_{i=1}^m$ be eigenvalues of the GUE such that 
$k_1 \sim n^{\gamma}$ where $0 < \gamma < 1$ and
$0 < k_{i+1} - k_i \sim n^{\theta_i}$, $0 < \theta_i < \gamma$.
Set
\begin{equation*}
X_i = \frac{x_{n-k_i} - \sqrt{2n} \left( 1 - 
\left( \frac{3 \pi k_i}{4 \sqrt{2} n} \right)^{2/3} \right)}{\left( 
\left( \frac{1}{12 \pi} \right)^{2/3} \frac{\log k_i}{n^{1/3} k_i^{2/3}}
\right)^{1/2}} \quad i=1, \ldots , m.
\end{equation*}  
Then as, $n \rightarrow \infty$, 
\begin{equation*}
\mathbb{P} \left[ X_1 \leq x_1, \ldots , X_m \leq x_m \right] 
\longrightarrow \Phi_{\Lambda}(x_1, \ldots ,x_m)
\end{equation*}
where $\Phi_{\Lambda}$ is the cdf for the 
$m$-dimensional normal distribution with covariance matrix
$\Lambda_{i,j}= 1 - \frac{1}{\gamma}\max \{\theta_k | \, i \leq k < j < m \}$
if $i<j$ and $\Lambda_{i,i} =1$.
\end{theorem}

\emph{Remark 1.} As one would expect the eigenvalues near the edge of the 
spectrum are less correlated than the ones in the bulk.\\

\emph{Remark 2.} The eigenvalues are quite correlated in the bulk. In order 
for $x_k$ and $x_m$ to be independent in the limit it must hold that
$|k-m| \sim n$. It is interesting to compare with the following
result by Mosteller\footnote{Mosteller actually allowed for $X_i$ to come 
from more general distributions.}, \cite{da} p. 201: \\
Let $X_i$, $i=1, \ldots, n$, be independent random variables  
uniformly distributed on $(0,1)$. Consider the asymptotic joint 
distribution of the $m$ sample quantiles $X_{n_j}$, $j=1,\ldots, m$,
where $n_j = [\lambda_j n] + 1$ and $0<\lambda_1 < \ldots <\lambda_m$.
\begin{theorem}[Mosteller]
As $n \rightarrow \infty$ the joint distribution of 
$X_{n_1}, \ldots, X_{n_m}$ tends to the 
$m$-dimensional normal distribution with means 
$\lambda_j$, variances $n^{-1}\lambda_j(1-\lambda_j)$ and correlations 
\begin{equation*}
\rho(X_{n_j}X_{n_{j'}}) = 
\sqrt{\frac{\lambda_j (1-\lambda_{j'})}{\lambda_{j'}(1-\lambda_j)}} 
\qquad j \leq j'.
\end{equation*}
\end{theorem} 
Hence, in this case $X_{n_1}, \ldots , X_{n_m}$ are globally correlated in 
the limit.

\section{Proofs of Theorems \ref{th1} and \ref{th2}}
\setcounter{equation}{0}

The proofs of Theorems \ref{th1} and \ref{th2} rely on a theorem due to
Costin, Lebowitz and Soshnikov, see \cite{cl} and \cite{so2}. 
Before presenting it we need some notation. 

Let $\{ P_t$\}, 
$t \in \mathbb{R_+}$, be a family
of random point fields, see \cite{so1}, on the real line such that their
correlation functions have a determinantal form\footnote{An example is the 
GUE.}. Call the determinant kernels $K_t(x,y)$ 
and let $\{ I_t \}$ be a set of intervals. We denote by $A_t$ the integral
operator on $L^2(I_t)$ defined by the kernel $K_t(x,y)$, 
$A_t : L^2(I_t) \rightarrow L^2(I_t)$. Furthermore, $\E_t$ and 
$\var_t$ 
denote the expectation and variance with respect to the probability measure
$P_t$. Finally, let $\#I_t$ stand for the number of particles in $I_t$.
\begin{theorem} [Costin-Lebowitz, Soshnikov] \label{cls}
Let $A_t = K_t \cdot \chi_{I_t}$ be a family of trace class operators
associated with the determinantal random point fields $\{P_t\}$ such 
that $\var_t(\#I_t) = \mathrm{Trace}(A_t - A_t^2)$ goes to infinity
as $t \rightarrow \infty$. Then 
\begin{equation*}
\frac{\#I_t - \E[\#I_t]}{\sqrt{\var(\#I_t)}} \longrightarrow N(0,1)
\end{equation*}
in distribution with respect to the random point field $P_t$.
\end{theorem}
The following lemmas will be proven in sections \ref{exp} and \ref{var}:
\begin{lemma} \label{exp1}
Let $t = t(k,n)$ be the solution to the equation
\begin{equation*}
n \frac{2}{\pi} \int_{-1}^{t} \sqrt{1-x^2} \, \udx = k .
\end{equation*}
where $k = k(n)$ is such that $k/n \rightarrow a \in (0,1)$
as $n \rightarrow \infty$.
The expected number of eigenvalues in the interval 
\begin{equation*}
I_n = \bigg[ \sqrt{2n} \, t + x \sqrt{\frac{\log n}{2n}}, \infty \bigg)
\end{equation*}
is given by
\begin{equation*}
\mathbb{E} \left[ \# I_n \right] = n - k -
\frac{x}{\pi} \sqrt{(1 - t^2) \log n} +
\mathcal{O} \left( \frac{\log n}{n} \right) .
\end{equation*}
\end{lemma}
\begin{lemma} \label{exp2}
The expected number of eigenvalues in the interval
$I_n =[ \sqrt{2n} \, t , \infty )$, where $t \rightarrow 1^-$ as
$n \rightarrow \infty$, is given by
\begin{equation*}
\mathbb{E} \left[ \# I_n \right] =
\frac{4 \sqrt{2}}{3 \pi} n(1-t)^{3/2} + \mathcal{O}(1) .
\end{equation*}
\end{lemma}
\begin{lemma} \label{var1}
Let $\delta > 0$ and suppose that $t$, which may 
depend on $n$, is such that $-1+ \delta \leq t < 1$ and
$n(1-t)^{3/2} \rightarrow \infty$ as $n \rightarrow \infty$.
Then the variance of the number of eigenvalues in the interval 
$I_n = [ t \sqrt{2n} , \infty )$ is given by 
\begin{equation*}
\var(\#I_n) = \frac{1}{2 \pi^2} \log [n(1-t)^{3/2}] \left( 1 +
\eta (n) \right)
\end{equation*}  
where $\lim_{n \rightarrow \infty} \eta (n) = 0 $.
\end{lemma}
Using these lemmas and Theorem \ref{cls} we are now ready to prove 
Theorems \ref{th1} and \ref{th2}.
\begin{proof}[Proof of Theorem \ref{th1}:]
Set
\begin{equation*}
I_n = \Bigg[ t \sqrt{2n} + 
\xi \left( \frac{\log n}{4(1-t^2)n} \right)^{1/2} , \infty \Bigg) .
\end{equation*}
Using Lemma \ref{exp1} and Lemma \ref{var1} we get
\begin{multline*}
\mathbb{P}_n \left[ 
\frac{x_k - t \sqrt{2n}}{\left( \frac{\log n}{4(1-t^2)n}
\right)^{1/2}} \leq \xi \right] = 
\mathbb{P}_n \left[ 
x_k \leq t \sqrt{2n} + \xi \left( \frac{\log n}{4(1-t^2)n}
\right)^{1/2} \right] \\ =
\mathbb{P}_n \left[ \# I_n \leq n-k \right] = 
\mathbb{P}_n \left[ \frac{\# I_n - \mathbb{E}_n [\# I_n]}
{\left( \var_n (\# I_n) \right)^{1/2}} \leq
\frac{n - k - \mathbb{E}_n [\# I_n]}{\left( \var_n (\# I_n)
\right)^{1/2}} \right] \\ =
\mathbb{P}_n \left[ \frac{\# I_n - \mathbb{E}_n [\# I_n]}
{\left( \var_n (\# I_n) \right)^{1/2}} \leq \xi + \epsilon (n) \right]
\end{multline*}
where $\epsilon (n) \rightarrow 0$ as $n \rightarrow \infty$.
By the Costin-Lebowitz-Soshnikov theorem the conclusion follows.
\end{proof}
\begin{proof}[Proof of Theorem \ref{th2}:]
Let $g(t)$ be the expected number of eigenvalues in the interval
$I_n = [t \sqrt{2n},\infty)$. We have that
\begin{multline*}
\mathbb{P}_n \left[ x_{n-k} \leq t \sqrt{2n} \right] =
\mathbb{P}_n \left[\# I_n \leq k \right] \\ =
\mathbb{P}_n \left[ \frac{\# I_n - g(t)}{(\var_n(\# I_n))^{1/2}}
\leq \frac{k - g(t)}{(\var_n(\# I_n))^{1/2}} \right].
\end{multline*}
If we can find $t$ such that 
\begin{equation} \label{limit}
\frac{k - g(t)}{(\var_n(\# I_n))^{1/2}} \rightarrow \xi
\end{equation}
as $n \rightarrow \infty$ then, by the Costin-Lebowitz-Soshnikov 
theorem, we are done.
The idea now is therefore to find a candidate for $t$. We will then 
insert this $t$ in the equation above to see if it is satisfied. 
Set for simplicity $h(t) = (\var_n(\# I_n))^{1/2}$. We get
from Lemma \ref{exp2} and Lemma \ref{var1} that 
\begin{eqnarray*}
& & g(t) = a_1 n(1-t)^{3/2} + \mathcal{O} (1) \\
& & h(t) = a_2 \log^{1/2} [n(1-t)^{3/2}] + o(\log^{1/2} [n(1-t)^{3/2}])
\end{eqnarray*}
where $a_i$ are
known constants. We need to study the equation
\begin{equation*}
k = g(t) + \xi h(t)
\end{equation*}
or, since $g$ is a strictly decreasing function,
\begin{equation*}
t = g^{-1}(k - \xi h(t)) \approx g^{-1}(k) - (g^{-1})'(k) \cdot \xi h(t) .
\end{equation*}
Since
\begin{equation*}
(g^{-1})'(k) = \frac{1}{g'(g^{-1}(k))}
\end{equation*}
we need to study $g^{-1}(k)$. 
\begin{eqnarray*} 
& & k \approx a_1 n(1-t)^{3/2} \Longrightarrow \\
& & t \approx 1 - \left(\frac{k}{a_1 n} \right)^{2/3}
\end{eqnarray*}
A reasonable guess for the derivative of $g$ is that
\begin{equation*}
g'(t) \approx - \frac{3 a_1}{2} n \sqrt{1-t}  .
\end{equation*}
We now get
\begin{equation*}
g'(g^{-1}(k)) \approx -\frac{3 a_1}{2} n \left( 
\left( \frac{k}{a_1 n} \right)^{2/3} \right)^{1/2} =
- \frac{3 a_1^{2/3}}{2} k^{1/3} n^{2/3} 
\end{equation*}
and 
\begin{equation*}
h(t) \approx h \left( g^{-1}(k) \right) \approx 
a_2 \log^{1/2} \left[ n \frac{k}{a_1 n} 
\right] \approx a_2 \log^{1/2} k
\end{equation*}
When gluing the pieces together one gets
\begin{equation*}
t \approx 1 - \left(\frac{k}{a_1 n} \right)^{2/3}
+ \xi \frac{2 a_2}{3 a_1^{2/3}} 
\frac{\log^{1/2} k}{k^{1/3} n ^{2/3}} .
\end{equation*}
When inserting this expression in (\ref{limit}) it turns out that it 
all works out. Some rearranging finally yields the result.
\end{proof}

\section{Proof of theorems \ref{th3} and \ref{th4}} \label{th34}
\setcounter{equation}{0}

We shall use the following theorem, \cite{so3}:
\begin{theorem} [Soshnikov]\label{multigaussian}
Let $(X,\mathcal{F},P_L)$ be a family of determinantal random point
fields with Hermitian locally trace class kernels $K_L$ and 
\begin{equation*}
\{I_L^{(1)}, \ldots, I_L^{(k)}\}_{L \geq 0}
\end{equation*}
be a family of Borel 
subsets of $\mathbb{R}$, disjoint for any fixed $L$, with compact 
closure. 

Suppose that the variance of the linear statistic 
$\sum_{i=-\infty}^{\infty} f_L(x_i)$, where 
\begin{equation*}
f_L(x) = \sum_{j=1}^k \alpha_j \cdot \chi_{I_L^{(j)}}(x) \qquad
\alpha_1, \ldots, \alpha_k \in \mathbb{R},
\end{equation*}
grows to infinity with $L$ in such a way that 
\begin{equation} \label{prereq}
\var_L \left( \#I_L^{(j)} \right) =  \mathcal{O} 
\left( \var_L \left(\sum_{i=-\infty}^{\infty} f_L(x_i) \right) \right)
\end{equation} 
for any $1 \leq j \leq k$. Then the Central limit theorem holds:
\begin{equation*}
\frac{\sum_{j=1}^{k} \alpha_j^{(L)} \#I_L^{(j)} -
\E_L \left[ \sum_{j=1}^{k} \alpha_j^{(L)} \#I_L^{(j)} \right]}
{\sqrt{\var_L \left( \sum_{j=1}^{k} \alpha_j^{(L)} \#I_L^{(j)} \right)}}
\longrightarrow N(0,1)
\end{equation*}
in distribution.
\end{theorem}
\emph{Remark 1:} The theorem in \cite{so3} is actually more general than the 
theorem stated here. \\

\emph{Remark 2:} If (\ref{prereq}) holds for any 
$\alpha_1, \ldots, \alpha_k$ then $\#I_L^{(1)}, \ldots,\#I_L^{(k)}$ are 
jointly normally distributed in the limit, \cite{du}. 
\begin{proof} [Proof of Theorem \ref{th3}:]
Take $k_i$, $s_i$, $\theta_i$ and $X_i$ as in the formulation of 
Theorem \ref{th3}. Note that $k_i - k_{i+1} \sim n^{\theta_i}$ implies that
$s_i - s_{i+1} \sim n^{\theta_i - 1}$. 

For any real numbers $x_i$ we have the identity (for $n$ large enough)
\begin{align*} 
\mathbb{P} [ X_1 \leq x_1, \ldots,  & X_m  \leq x_m ]  \\  
= \mathbb{P} \Bigg[ &
\frac{\# I_1 - \E [\# I_1]}{(\var(\# I_1))^{1/2}} \leq 
\frac{n-k_1-\E [\# I_1]}{(\var(\# I_1))^{1/2}}, \\
& \frac{\# I_1 + \# I_2 - 
\E [\# I_1 + \# I_2]}{(\var(\# I_1 + \# I_2))^{1/2}} 
\leq
\frac{n-k_2-\E [\# I_1 + \# I_2]}{(\var(\# I_1 + 
\# I_2))^{1/2}}, 
\ldots , 
\end{align*}
\begin{align*}
& \frac{ \sum_{i=1}^m \# I_i - 
\E [\sum_{i=1}^m \# I_i]}
{\left( \var \left( \sum_{i=1}^m \# I_i \right) \right)^{1/2}} 
\leq
\frac{n-k_m-\E [\sum_{i=1}^m \# I_i]}
{\left( \var \left( \sum_{i=1}^m \# I_i \right) \right)^{1/2}} \Bigg] .
\end{align*}
Here the intervals $I_i$ are given by
\begin{align*}
I_1 & = \left( s_1 \sqrt{2n} + 
x_1 \left( \frac{\log n}{4(1-s_1^2)n} \right)^{1/2}, \infty \right) \\
I_i & = 
\Bigg(s_i \sqrt{2n} + x_i \left( \frac{\log n}{4(1-s_i^2)n} \right)^{1/2}, \\
& \, \, \, \, \, \, \, \, \, \, \, \, \, \, \, \, \, \, \, \, \, \, \, \,
\, \, \, \, \, \, \, \, \, \, \, \, \, \, \, \, \, \, \, \, \, \, \, \,
\, \, \, \, \, \, \, \, \, \, \, \, \, \, \, \, \, \, \, \, \, \, \, \,
\, \, \, \, \, \, \, \, \, \, \, \, \, 
s_{i-1} \sqrt{2n} + x_{i-1} \left( \frac{\log n}{4(1-s_{i-1}^2)n} 
\right)^{1/2} \Bigg] 
\end{align*}
where $2 \leq i \leq m$. We will now show that the random variables
\begin{equation*}
\#I_1 , \#I_1 + \#I_2, \ldots , \sum_{i=1}^m \#I_i 
\end{equation*}
are jointly normal in the limit. To this end we shall use Theorem 
\ref{multigaussian} above to show that all linear 
combinations of the variables are normally distributed in the limit. Since 
\begin{equation*}
\alpha_1 \#I_1 + \alpha_2 (\#I_1 + \#I_2) =
(\alpha_1 + \alpha_2) \#I_1 + \alpha_2 \#I_2
\end{equation*}
and so forth
it is clear that one can instead look at all linear combinations of 
$\{ \#I_i \}_1^m$.
Hence, by Theorem \ref{multigaussian}\footnote{The theorem by Soshnikov 
does not apply directly to this situation since $I_1$ does not have 
compact closure. This is however easily overcome simply by chopping of 
the interval far out where the probability of finding any eigenvalue
is exponentially small in $n$.} we must calculate (see Appendix B)
\begin{multline*}
\var \left( \alpha_1 \#I_1 + \alpha_2 \#I_2 + \ldots +
\alpha_m \#I_m \right) \\
= \sum_{i=1}^m \alpha_i^2 
\iint_{I_i \times I_i^c} K_n^2(x,y) \, \udx \udy  \\
 - \sum_{i \neq j}^m \alpha_i \alpha_j 
\iint_{I_i \times I_j} K_n^2(x,y) \, \udx \udy
\end{multline*}
to see that it is of magnitude $\log n$. 
First define the set $M$ by $k \in M \iff \theta_k = 1$. Hence
\begin{equation*}
M = \{ k_1, \ldots , k_j\} ; \quad
1 \leq k_1 < k_2 < \ldots < k_j \leq m-1 
\end{equation*}
for some $j$ such that $0 \leq j \leq m-1$. 
Suppose first that $M$ is empty which means that $\theta_i < 1$ for all $i$.
If $\alpha_1 \neq 0$ then, by using the inequality 
\begin{equation*}
x y \leq \frac{1}{2}(x^2 + y^2), 
\end{equation*}
we get
\begin{eqnarray} \nonumber
& & \var \left( \alpha_1 \#I_1 + \alpha_2 \#I_2 + \ldots +
\alpha_m \#I_m \right) \\ \nonumber
& & \qquad \qquad \geq \sum_{i=1}^m \alpha_i^2 
\iint_{I_i \times I_i^c} K_n^2(x,y) \, \udx \udy  \\ \nonumber
& & \qquad \qquad \qquad \qquad \qquad \qquad - 
\sum_{i \neq j}^m \frac{1}{2}(\alpha_i^2 + \alpha_j^2) 
\iint_{I_i \times I_j} K_n^2(x,y) \, \udx \udy  \\ \label{vars}
& & \qquad \qquad \qquad = \sum_{i=1}^m \alpha_i^2 \bigg(
\iint_{I_i \times I_i^c} K_n^2(x,y) \, \udx \udy  \\ \nonumber
& & \qquad \qquad \qquad \qquad \qquad \qquad \qquad -\sum_{j \neq i}^m
\iint_{I_i \times I_j} K_n^2(x,y) \, \udx \udy \bigg).
\end{eqnarray} 
All the terms under the $i$-summation in (\ref{vars}) are non-negative 
and the first term can be 
calculated as in the proof of Lemma \ref{var1} which we provide in section
\ref{var}.
In the proof of this lemma we show that in the domain
\begin{equation*}
\Omega = 
\left\{ (x,y) ; s \leq x \leq s + \frac{1}{\log n} , 
s - \frac{1}{\log n} \leq y \leq s \right\}
\end{equation*}
it holds that 
\begin{equation*}
2n K_n(\sqrt{2n} x,\sqrt{2n} y) = \frac{1}{2 \pi^2 (x-y)^2} +
\mathcal{O} \left( \frac{1}{\log n} \right)
\end{equation*}
It is also shown that if 
\begin{equation*}
\Omega' = \left\{ (x,y) ; \sqrt{2n} s \leq x \leq \infty, - \infty < y \leq
\sqrt{2n} s \right\} \not \quad \sqrt{2n} \cdot \Omega
\end{equation*}
then\footnote{Since $K_n(x,y) = K_n(y,x)$ it is clear that the same
estimates hold in the domains obtained from reflection with respect to
the $x=y$-line.}
\begin{equation*}
\iint_{\Omega'} K_n^2(x,y) \, \udx \udy = \mathcal{O} (\log \log n).
\end{equation*}
In what follows we shall often make use of these facts without 
mentioning it. The main contribution to the first term in (\ref{vars})
can now be calculated to be (disregarding $\alpha_1^2$)
\begin{equation*}
\int_{s_1}^{s_1 + \frac{1}{\log n}} 
\int_{s_1 - \frac{1}{\log n}}^{s_1 - n^{\theta^* - 1}} \frac{1}{(x-y)^2}
\, \udy \udx = \frac{1-\theta^*}{2 \pi^2} \log n + \mathcal{O} (\log \log n)
\end{equation*}
where $\theta^* = \max_{i} \theta_i < 1$. By our definition of $\sim$ 
above the integration in the $y$-variable should have been 
over the interval $\left(s_1 - 1/\log n, s_1 - h(n) n^{\theta^* -1} \right)$ where
$h(n)$ satisfies (\ref{sim}). However, because of the logarithmic 
answer this $h$ will only produce lower order terms.

Now suppose that $j=0$ as before, $\alpha_1 = \ldots = \alpha_{k-1} = 0$
but $\alpha_k \neq 0$. In this case we get
\begin{multline*}
\var \left( \alpha_k \#I_k + \ldots + \alpha_m \#I_m \right) \\
\geq \sum_{i=k}^m \alpha_i^2 \bigg(
\iint_{I_i \times I_i^c} K_n^2(x,y) \, \udx \udy \\
- \sum_{k \leq j \neq i}^m
\iint_{I_i \times I_j} K_n^2(x,y) \, \udx \udy \bigg).
\end{multline*} 
Using the estimates above it is straightforward to verify that the
$k$-term is of order $\log n$.

When $j \geq 1$ meaning that there is at least
one $k$ with $\theta_k = 1$, things are only slightly more complicated. 
Let $k^*$ be the largest integer $i$ such that $\theta_i = 1$. It is 
sufficient to consider the case when there exists $i \geq k^* + 1$
such that $\alpha_i \neq 0$. On the other hand if this is the case then
we are in a situation very similar to when $j=0$. Either 
$\alpha_{k^*+1} \neq 0$ or $\alpha_{k^*+1} = \ldots = \alpha_{l-1} = 0$
but $\alpha_{l} \neq 0$. The details are left out.

It is hence a fact that 
\begin{equation*}
\#I_1 , \#I_1 + \#I_2, \ldots , \sum_{i=1}^m \#I_i 
\end{equation*} 
in the limit have a joint normal distribution.

To complete the proof we need to calculate the correlations between the
different $\#I_i$'s. 
If $j < i$ we have that $s_j - s_i \sim n^{-\gamma}$
where $\gamma = 1 - \max_{j \leq k < i} \theta_k$. Set
\begin{equation*}
X_k = \sum_{m=1}^k \#I_m . 
\end{equation*}
From a straightforward calculation (as above) we get that
\begin{multline*}
\var (X_i - X_j)  = \var\left( \sum_{k=j+1}^i \#I_k \right)  
= \var \left( \# \bigcup_{k=j+1}^i I_k \right) \\ = 
\frac{1-\gamma}{\pi^2} \log n + \mathcal{O} (\log \log n).
\end{multline*}
Since
\begin{equation*}
\var(X_k) = \frac{1}{2 \pi^2} \log n + \mathcal{O} (\log \log n)
\end{equation*}
the correlation $\rho$ is given by
\begin{equation*}
\rho(X_i,X_j) = \frac{\frac{1}{2} \left( \var(X_j) + \var(X_i) -
\var(X_i - X_j) \right)}{ \sqrt{\var(X_i)\var(X_j)}} =
\gamma + o (1).
\end{equation*}
\end{proof}

\begin{proof}[Proof of Theorem \ref{th4}.]
This proof is of course very similar to the previous one so some
details will be skipped. 

With notation as in the formulation of Theorem \ref{th4} 
the intervals of interest (confer with the previous proof) are in this case 
\begin{eqnarray*}
& & I_1 = \left( \sqrt{2n} \left( 1 - C_1 \left(\frac{k_1}{n} \right)^{2/3}
\right) + x_1 C_2 
\left( \frac{\log k_1}{n^{1/3} k_1^{2/3}} \right)^{1/2}, \infty \right) \\
& & I_i = \Bigg(
\sqrt{2n} \left( 1 - C_1 \left(\frac{k_{i}}{n} \right)^{2/3} \right) + 
x_i C_2 \left( \frac{\log k_i}{n^{1/3} k_i^{2/3}} \right)^{1/2}, \\
& & \qquad \qquad \qquad \sqrt{2n} 
\left( 1 - C_1 \left(\frac{k_{i-1}}{n} \right)^{2/3} \right) + x_{i-1} C_2 
\left( \frac{\log k_{i-1}}{n^{1/3} k_{i-1}^{2/3}} \right)^{1/2} \Bigg]
\end{eqnarray*}
where $C_1$, $C_2$ are known constants and $2 \leq i \leq m$. 
Given any $x_1, \ldots, x_m$ it is easy to see
that the sets $I_i$ are intervals if $n$ is sufficiently large.
As in the previous proof we want to show that 
\begin{equation*}
\#I_1 , \#I_2, \ldots , \#I_m 
\end{equation*}
are jointly normally distributed. The way to prove this is the same 
as before but some details are different. By Lemma \ref{var1} we need to 
show that
\begin{equation*}
\log n = \mathcal{O} 
\left(\var\left(\sum_{i=1}^m \alpha_i \#I_i \right) \right) 
\end{equation*} 
for any real $\alpha_i$'s such that $\alpha_i \neq 0$ for some $i$. 

Let $t=t(n)$ be such that
$n^{\epsilon - \frac{2}{3}} \leq 1-t \leq n^{- \epsilon}$ for some
$0 < \epsilon \leq 1/3 $.
From the proof of Lemma \ref{var1} below one can see that in the set
\begin{equation*}
\Omega_t = \left\{ (x,y) ; \, t \leq x \leq t + \frac{1-t}{\log n},
t - \frac{1-t}{\log n} \leq y \leq t \right\}
\end{equation*}
it holds that 
\begin{equation*}
2n K_n^2(\sqrt{2n}x, \sqrt{2n}y) = \frac{1}{2 \pi^2 (x-y)^2} + 
\mathcal{O} \left( \frac{1}{\log n} \right) .
\end{equation*}
Returning to the variance calculation we first assume that $\alpha_1 \neq 0$.
We know from the previous proof that in this case it is sufficient to
to show that
\begin{equation*}
\iint_{I_1 \times \left(I_1^c \setminus \bigcup_{i=2}^m I_i \right)} K_n^2(x,y) 
\, \udy \udx
\end{equation*}
is of order $\log n$. In fact since the integrand is non-negative
it is enough if
\begin{equation*}
\iint_{I^* \times I_*} \frac{1}{(x-y)^2} \, \udy \udx
\end{equation*}
is of order $\log n$ where 
\begin{align*}
I^* &= \left( t_1 + r_1, t_1 + \frac{1-t_1}{\log n} \right) \\
I_* &= \left( t_1 - \frac{1-t_1}{\log n}, t_1 \right)
\end{align*}
and 
\begin{align*}
t_i &= 1 - C_1 \left(\frac{k_{i}}{n} \right)^{2/3} \\
r_i &= x_i C_2 \left( \frac{\log k_i}{n^{1/3} k_i^{2/3}} \right) .
\end{align*}
An elementary calculation shows that this integral is indeed of 
order $\log n$. 

If $\alpha_1 = \ldots = \alpha_{k-1} = 0$ but $\alpha_k \neq 0$
it is sufficient that the integral
\begin{equation*}
\iint_{J^* \times J_*} \frac{1}{(x-y)^2} \, \udy \udx
\end{equation*}
is of order $\log n$ where
\begin{align*}
J^* &= \left( t_{k-1} + r_{k-1}, 
t_{k-1} + \frac{1-t_{k-1}}{\log n} \right) \\
J_* &= \left( t_k , t_{k-1} \right).
\end{align*}
Again we get the size $\log n$.
This proves that we get a normal distribution in the limit. The 
calculations of the correlations are very similar to
the bulk case and the details are not presented here. 
\end{proof}

\section{The expected number of eigenvalues in $I_n$} \label{exp}
\setcounter{equation}{0}

In this section and the next we shall need  
asymptotics for the Airy function and the Hermite polynomials.
In \cite{dkmvz} the asymptotics for a class containing the Hermite case was
studied. It is shown that there exists a $\delta > 0$ such that the following
holds: 

\noindent $\textrm{1. } -1 + \delta \leq x \leq 1 - \delta$
\begin{align*}
h_n  \left(  \sqrt{2n}  x \right) & e^{-n x^2}  \\
& = \sqrt{ \frac{2}{\pi \sqrt{2n}}} \frac{1}{(1-x^2)^{1/4}}
\left( \cos \left[ 2n F(x) - \frac{1}{2} \arcsin (x) \right] +
\mathcal{O} (n^{-1}) \right) 
\end{align*}
\noindent $\textrm{2. } 1  -  \delta \leq x < 1$ 
\begin{align*}
h_n  \left( \sqrt{2n} x \right)&  e^{-n x^2} 
= (2n)^{-1/4} \\ \times \bigg\{ & 
\left( \frac{1+x}{1-x} \right)^{1/4} \left[ 3nF(x) \right]^{1/6}
\ai \left( -\left[ 3nF(x) \right]^{2/3} \right)
(1 + \mathcal{O}(n^{-1})) \\
- & \left( \frac{1-x}{1+x} \right)^{1/4} \left[ 3nF(x) \right]^{-1/6}
\dai \left( -\left[ 3nF(x) \right]^{2/3} \right)
(1 + \mathcal{O}(n^{-1})) \bigg\} 
\end{align*}
\noindent $\textrm{3. } 1 < x \leq 1 + \delta$
\begin{align*}
h_n(\sqrt{2n} x) & e^{-n x^2} 
= (2n)^{-1/4} \\ \times \bigg\{ &
\left( \frac{x+1}{x-1} \right)^{1/4} \left[ 3nF(x) \right]^{1/6}
\ai \left( \left[ 3nF(x) \right]^{2/3} \right)  \\
- &\left( \frac{x-1}{x+1} \right)^{1/4} \left[ 3nF(x) \right]^{-1/6}
\dai \left( \left[ 3nF(x) \right]^{2/3} \right)  \bigg\}
(1 + \mathcal{O}(n^{-1})) 
\end{align*}
$\textrm{4. } x > 1 + \delta$
\begin{equation*}
h_n(\sqrt{2n} x) e^{-n x^2} =
\mathcal{O} \left( n^{-1/4} e^{- n F(x)} \right)
\end{equation*}
In these expressions $\ai$ stands for the Airy function and 
\begin{equation} \label{F}
F(x) = \left| \int_{x}^{1} \sqrt{|1-y^2|} \, \udy \right|.
\end{equation}
There are of course also similar asymptotics for the Hermite polynomials
near the point $-1$. 

The Airy function is bounded on the real line. It is exponentially small
in $x$ on $\mathbb{R}_+$  and for $r>0$ it holds that, \cite{ol},
\begin{align*}
\ai(-r) & = \pi^{-1/2} r^{-1/4} \bigg\{
\cos \left[ \frac{2}{3} r^{3/2} - \frac{\pi}{4} \right] + 
\mathcal{O} (r^{-3/2}) \bigg\} \\
\dai(-r) & = \pi^{-1/2} r^{1/4} \bigg\{
\sin \left[ \frac{2}{3} r^{3/2} - \frac{\pi}{4} \right] + 
\mathcal{O} (r^{-3/2}) \bigg\} .
\end{align*}
\begin{proof}[Proof of Lemma \ref{exp1}:]
Set
\begin{equation*}
f_n(t) = t + x \frac{\sqrt{\log n}}{2n}.
\end{equation*}
We have that
\begin{equation*}
\mathbb{E} \left[ \# I_n \right] =
\int_{f_n(t)}^{\infty} n \rho_n(x) \, \udx
\end{equation*}
where $\rho_n$ is the scaled density for the eigenvalues
(the limiting density has support in $[-1,1]$).
From symmetry one gets 
\begin{equation*}
\int_{f_n(t)}^{\infty} n \rho_n(x) \, \udx =
\frac{n}{2} - \int_0^{f_n(t)} n \rho_n(x) \, \udx .
\end{equation*}
Formula (4.2) in \cite{em} applied to the Hermitian case says that for
$\delta > 0$ it holds that
\begin{align*}
n \rho_n(x) = & n \frac{2}{\pi} \sqrt{1-x^2} \\
& + \frac{1}{4 \pi} \left( \frac{1}{x-1} - \frac{1}{x+1} \right)
\cos \left[ n \frac{2}{\pi} \int_x^1 \sqrt{1-y^2} \, \udy \right]
+ \mathcal{O} (n^{-1}) 
\end{align*} 
if $x \in [-1 + \delta, 1 - \delta]$. We get that
\begin{multline*}
\mathbb{E} \left[ \# I_n \right] =
\frac{n}{2} - n \frac{2}{\pi} \int_0^{f_n(t)} \sqrt{1-x^2} \, \udx +
\mathcal{O} (n^{-1}) \\ = 
 n - n \frac{2}{\pi} \int_{-1}^{f_n(t)} \sqrt{1-x^2} \, \udx + 
\mathcal{O} (n^{-1}) \\ =
 n - n \frac{2}{\pi} \left( \int_{-1}^{t} \sqrt{1-x^2} \, \udx +
\sqrt{1-t^2} x \frac{\sqrt{\log n}}{2n} + 
\mathcal{O} \left( \frac{\log n}{n^2} \right) \right) +
\mathcal{O} (n^{-1}) \\ =
n - k - \frac{x}{\pi} \sqrt{(1-t^2) \log n} +
\mathcal{O} \left( \frac{\log n}{n} \right) .
\end{multline*} 
\end{proof}
\begin{proof} [Proof of Lemma \ref{exp2}:]
From the formulas (4.4) and (4.21) 
in \cite{em} one gets after some minor calculations that
\begin{multline*}
n \rho_{n}(x) = \left( \frac{\Phi'(x)}{4 \Phi(x)} -
\frac{\gamma'(x)}{\gamma(x)} \right) 
[ 2 \ai(\Phi(x)) \dai(\Phi(x))] \\
+ \Phi'(x) \left[ ( \dai(\Phi(x)))^2 - \Phi(x) (\ai(\Phi(x)))^2 
\right] + \mathcal{O} \left(\frac{1}{n (\sqrt{1-x})}\right) 
\end{multline*} 
in a fixed neighborhood of $[0,1]$.
As before $\rho_n$ is the scaled density for the eigenvalues so that 
\begin{equation*}
g(t) = \int_{t}^{\infty} n \rho_n(x) \, \udx 
\end{equation*}
and the functions $\gamma$ and $\Phi$ are given by
\begin{eqnarray*}
\gamma(x) & = & \left( \frac{x-1}{x+1} \right)^{1/4} \\
\Phi(x) & = & \left\{ \begin{array}{ll}
- \left( 3n \int_{x}^{1} \sqrt{1-y^2} \, \udy \right)^{2/3} & 
\textrm{if $x \leq 1$} \\
\quad \left( 3n \int_{1}^{x} \sqrt{y^2-1} \, \udy \right)^{2/3} & 
\textrm{if $x > 1$.} 
\end{array} \right.
\end{eqnarray*}
The function $\gamma$ is evaluated by taking the limit from the upper half 
plane using the principal branch. 

The fact that the asymptotics only holds
for $x \in [0,1 + \delta]$ for some $\delta > 0$ (independent of $n$) is not
a problem. It is not difficult to show that, for $x \geq 1+\delta$, 
$\rho_n(x)$ is exponentially small in $n$ and exponentially decaying
in $x$. 

We shall now study the different terms in the asymptotical expression for 
$\rho_n$ above. When looking at the asymptotics for $\ai$ and $\dai$ it
easy to see that 
\begin{equation*}
\left| \ai(x) \dai(x) \right| = \mathcal{O} (1) .
\end{equation*}
A calculation shows that
\begin{equation*}
\left( \frac{\Phi'(x)}{4 \Phi(x)} -
\frac{\gamma'(x)}{\gamma(x)} \right) =
\mathcal{O} (1)
\end{equation*}
hence
\begin{equation*}
\int_{t}^{1+\delta} \left( \frac{\Phi'(x)}{4 \Phi(x)} -
\frac{\gamma'(x)}{\gamma(x)} \right) 
[ 2 \ai(\Phi(x)) \dai(\Phi(x))] \, \udx = \mathcal{O} (1).
\end{equation*}
The main contribution comes from the second term. In fact a primitive 
function can be found for this expression (see Appendix A):
\begin{align*}
\int_{t}^{1+\delta} 
\Phi'(x) & \left[ \dai^2(\Phi(x)) - \Phi(x) \ai^2(\Phi(x)) 
\right] \, \udx \\ & = \big[\mathrm{set} \, \, y = \Phi(x) \big]   \\
&= \int_{\Phi(t)}^{\Phi(1+\delta)} \dai^2(y) - y \ai^2(y) \, \udy \\
& = - \left[ \frac{2}{3} \left( y^2 \ai^2(y) - y \dai^2(y) \right)
- \frac{1}{3} \ai(y) \dai(y) \right]_{\Phi(t)}^{\Phi(1+\delta)} \\
&= \frac{2}{3} \left( \Phi^2(t) \ai^2(\Phi(t)) -  \Phi(t)
\dai^2(\Phi(t)) \right)  \\ 
& \qquad \qquad \qquad \qquad \qquad 
- \frac{1}{3} \ai(\Phi(t))\dai(\Phi(t)) 
+ \mathcal{O} \left( e^{-c n} \right)
\end{align*}
Here $c$ is a positive constant. Integrating the third term only
gives a contribution of order $n^{-1}$. One can now use the 
asymptotics for the Airy function and its derivative to get
the stated result.
\end{proof}

\section{The variance of the number of eigenvalues in $I_n$} \label{var}
\setcounter{equation}{0}

\begin{proof}[Proof of Lemma \ref{var1}:]
The proof will be divided into two basic cases. The first case is 
when $1-t > \delta$ for a fixed $\delta > 0$, that is, we are in the bulk.
The second case is when $t = t(n) \rightarrow 1^-$ as 
$n \rightarrow \infty$, that is, we are close to the right spectrum edge.

First define $I_n = [ t \sqrt{2n} , \infty )$ and $\# I_n$ as the
number of eigenvalues in $I_n$. It is a fact (see Appendix B) that
\begin{equation*}
\var(I_n) = \int_{I_n}{ \int_{I_n^c}{K_n^2(x,y)dx} dy}
\end{equation*}
where $K_n$ is the Hermite kernel.
The advantage with this representation is that there is only one
singular point involved in the Christoffel-Darboux representation of $K_n(x,y)$:
\begin{equation*}
K_n(x,y) = \sqrt{\frac{n}{2}} 
\frac{h_n(x) h_{n-1}(y) - h_{n-1}(x) h_n(y) }{x - y}
e^{-\frac{1}{2} (x^2 + y^2)}
\end{equation*}
\textbf{Case I (the bulk)}:
After a change of variables ($x \rightarrow \sqrt{2n} x$) we
get the integrand
\begin{equation*}
\left[ \sqrt{2n} \, K_n (\sqrt{2n} x, \sqrt{2n} y) \right]^2 .
\end{equation*}
First consider the domain where both variables are in the Bulk:
\begin{equation} \label{Gamma}
\Gamma = \left\{ (x,y) ; t \leq x \leq 1 - \delta , 
-1 + \delta \leq y \leq t \right\}
\end{equation} 
Recall that in $\Gamma$, $h_n$ has asymptotics given by
\begin{multline*}
h_n ( \sqrt{2n} x ) e^{-n x^2} \\ =
\left( \frac{2}{\pi \sqrt{2n}} \right)^{1/2} \frac{1}{(1-x^2)^{1/4}}
\left( \cos \left[ 2n F(x) - \frac{1}{2} \arcsin (x) \right] +
\mathcal{O} (n^{-1}) \right) 
\end{multline*}
where
\begin{equation*}
F(x) = \int_x^1 \sqrt{1-z^2} \, \udz =
\frac{1}{2} \left( \arccos x - x \sqrt{1-x^2} \right) .
\end{equation*}
The asymptotics for $h_{n-1}$ becomes 
\begin{align*}
h_{n-1} & ( \sqrt{2n} x ) e^{-n x^2} =
\left( \frac{2}{\pi \sqrt{2(n-1)}} \right)^{1/2} \\ & \qquad \times 
\frac{1}{(1-x_n^2)^{1/4}} 
\left( \cos \left[ 2(n-1) F(x_n) - \frac{1}{2} \arcsin (x_n) \right] +
\mathcal{O} (n^{-1}) \right) \\ & =
\left( \frac{2}{\pi \sqrt{2n}} \right)^{1/2}
\frac{1}{(1-x^2)^{1/4}} \\ & \qquad \qquad \qquad \qquad
\times \left( \cos \left[ 2(n-1) F(x_n) - \frac{1}{2} \arcsin (x) \right] +
\mathcal{O} (n^{-1}) \right)
\end{align*}
where $x_n = \sqrt{\frac{n}{n-1}} x$.
A Taylor expansion gives
\begin{equation*}
F (x_n) = F(x) - \frac{x}{2(n-1)} \sqrt{1-x^2} + \mathcal{O} (n^{-2})
\end{equation*}
and hence 
\begin{multline*}
2 (n-1) F(x_n) 
= 2n F(x) - 2F(x) -x \sqrt{1-x^2} + \mathcal{O} (n^{-1}) \\ 
= 2n F(x) - \arccos x + \mathcal{O}(n^{-1}) .
\end{multline*}
One can now write
\begin{align*}
h_n & (\sqrt{2n} x ) h_{n-1} ( \sqrt{2n}y )
e^{-n (x^2+y^2)} \\ & =
\frac{2}{\pi \sqrt{2n} (1-x^2)^{1/4} (1-y^2)^{1/4}}  \\ & \qquad \times
\cos \left[ 2n F(x) - \frac{1}{2} \arcsin x \right]
\cos \left[ 2n F(y) - \frac{1}{2} \arcsin y - \arccos y \right]
\\ & \qquad \qquad \qquad \qquad \qquad \qquad \qquad \qquad \qquad \qquad
\qquad \qquad \qquad
+ \mathcal{O} (n^{-3/2}) .
\end{align*}
Set
\begin{align*}
\alpha_x & = 2n F(x) - \frac{1}{2} \arcsin x \\
\theta_x & = \arccos x .
\end{align*}
By the Christoffel-Darboux formula 
\begin{multline*}
\sqrt{2n} \, K_n (\sqrt{2n}x,\sqrt{2n}y)  =
\frac{1}{\pi (1-x^2)^{1/4} (1-y^2)^{1/4}}  \\
\times \frac{\cos \alpha_x \cos [\alpha_y - \theta_y] -
\cos [\alpha_x - \theta_x] \cos \alpha_y + \mathcal{O}(n^{-1})}{x-y}.
\end{multline*} 
To prepare for integration we now divide $\Gamma$ into four disjoint 
sets. Set
\begin{eqnarray*}
& & \Gamma_0 = \left\{ (x,y) ; t \leq x \leq t + \frac{1}{n} ,
t - \frac{1}{n} \leq y \leq t \right\} \\ 
& & \Gamma_1 = \Gamma_1^1 \cup \Gamma_1^2 =
\left\{ (x,y) ; t \leq x \leq t + \frac{1-t}{r(n)} ,
t - \frac{t+1}{r(n)} \leq y \leq t - \frac{1}{n} \right\} \bigcup \\
& & \qquad \qquad \qquad \quad
\left\{ (x,y) ; t + \frac{1}{n} \leq x \leq t + \frac{1-t}{r(n)} ,
t - \frac{1}{n} \leq y \leq t \right\} \\
& & \Gamma_2 = \Gamma \setminus \left( \Gamma_0 \cup \Gamma_1 \right)
\end{eqnarray*}
where $r(n) = \log n$ and $\Gamma$ was defined in (\ref{Gamma}).

$\Gamma_0$: When integrating over $\Gamma_0$ one can use the fact that
\begin{equation*}
\sqrt{2n} \, |K_n (\sqrt{2n} x, \sqrt{2n} y)| \leq
C n \, \left| \frac{\sin(x-y)}{x-y} \right| 
\end{equation*}
where $C>0$. Hence
\begin{equation*}
\int_{\Gamma_0} \left[ \sqrt{2n} \, K_n(\sqrt{2n} x, \sqrt{2n} y)
\right]^2 \, \udx \udy = \mathcal{O}(1) .
\end{equation*}
\indent $\Gamma_1$: For $x \in \Gamma_1$ it holds that
\begin{equation*}
\theta_x = \arccos x = \arccos t + \mathcal{O} 
\left( \frac{1}{r(n)} \right) 
\end{equation*}
and of course also the equivalent for $\theta_y$. Defining
$\theta = \arccos t$ we get by the use of some trigonometric identities
that
\begin{align*}
\cos \alpha_x \cos [\alpha_y - \theta_y] & -
\cos [\alpha_x - \theta_x] \cos \alpha_y \\ = &
\cos \alpha_x \cos [\alpha_y - \theta] -
\cos [\alpha_x - \theta] \cos \alpha_y  + 
\mathcal{O}\left(\frac{1}{r(n)}\right) \\ = &
\sqrt{1-t^2} \sin [ \alpha_y - \alpha_x ] + 
\mathcal{O} \left( \frac{1}{r(n)} \right).
\end{align*}
Since 
\begin{equation*}
\frac{\sqrt{1-t^2}}{(1-x^2)^{1/4} (1-y^2)^{1/4}} = 
1 + \mathcal{O} \left( \frac{1}{r(n)} \right)
\end{equation*}
and
\begin{equation*}
\alpha_x - \alpha_y = 2n (F(x) - F(y)) + 
\mathcal{O} \left( \frac{1}{r(n)}\right)
\end{equation*}
we now get
\begin{align*}
\iint_{\Gamma_1} & \left[ \sqrt{2n} \, 
K_n(\sqrt{2n} x, \sqrt{2n} y) \right]^2 \, \udx \udy \\ = &
\iint_{\Gamma_1^1} \frac{1}{\pi^2} 
\frac{\sin^2[2n(F(x)-F(y))] + \mathcal{O} \left( \frac{1}{r(n)} 
\right)}{(x-y)^2} \, \udx \udy +
\iint_{\Gamma_1^2} \frac{\mathcal{O}(1)}{(x-y)^2}\, \udx \udy \\ = &
\frac{1}{2 \pi^2} \iint_{\Gamma_1^1} 
\frac{1 - \cos \left[ 4n(F(x)-F(y))\right]}{(x-y)^2} \, \udx \udy +
\mathcal{O} \left( \log r(n) \right) \\ = &
\frac{1}{2 \pi^2} \log n - \frac{1}{2 \pi^2} \iint_{\Gamma_1^1} 
\frac{\cos \left[ 4n(F(x)-F(y))\right]}{(x-y)^2} \, \udx \udy  +
\mathcal{O} \left( \log r(n) \right) .
\end{align*}
The remaining integral is not bigger than a constant as will now be 
shown. A partial integration in the y-variable gives
\begin{multline*}
\iint_{\Gamma_1^1} 
\frac{\cos \left[ 4n(F(x)-F(y))\right]}{(x-y)^2} \, \udx \udy \\ =
\int_{t}^{t+\frac{1-t}{r(n)}} \Bigg(
\left[ \frac{\sin \left[ 4n(F(x)-F(y)) \right]}{4nF'(y) (x-y)^2} 
\right]_{t-\frac{t+1}{r(n)}}^{t-\frac{1}{n}} \\ -
\int_{t-\frac{t+1}{r(n)}}^{t-\frac{1}{n}} 
\sin \left[ 4n(F(x)-F(y)) \right]
\left( \frac{1}{4n \left[ F'(y)(x-y)^2 \right]} \right)_y' 
\, \udy \Bigg) \, \udx \\ 
=: I_1 - I_2 .
\end{multline*}
Both these integrals are easy to estimate:
\begin{eqnarray*}
|I_1| \leq C \int_{t}^{t+\frac{1-t}{r(n)}} 
\frac{1}{n(x-(t-n^{-1}))^2} \, \udx = \mathcal{O}(1).
\end{eqnarray*}
We have
\begin{equation*}
\left( \left[ F'(y)(x-y)^2 \right]^{-1} \right)_y' =
- \frac{y}{(1-y^2)^{3/2} (x-y)^2} - \frac{2}{\sqrt{1-y^2} (x-y)^3}
\end{equation*}
which gives
\begin{equation*}
|I_2| \leq C \iint_{\Gamma_1^1} \frac{1}{n(x-y)^3} =
\mathcal{O} \left( 1 \right) .
\end{equation*}
$\Gamma_2$: In $\Gamma_2$ it holds that 
\begin{equation*}
\left[ \sqrt{2n} \, K_n(\sqrt{2nx},\sqrt{2ny}) \right]^2 =
\mathcal{O} \left( \frac{1}{(x-y)^2} \right)
\end{equation*}
and a trivial calculation gives
\begin{equation*}
\iint_{\Gamma_2} \frac{1}{(x-y)^2} \, \udx \udy  =
\mathcal{O} \left( \log r(n) \right) .
\end{equation*} 
To complete case I we must also integrate over
$I_n \times I_n^c \setminus \Gamma$. By using the appropriate asymptotics
for the Hermite polynomials, the Airy function and its derivative and
then taking absolute values in the integrals it is straight forward to
show that the contribution from this domain is $\mathcal{O}(1)$. We omit the
details.

\textbf{Case II (the spectrum edge)}: First consider the subdomain 
\begin{equation*}
\Omega = \left\{ (x,y); t \leq x \leq 1-C n^{-1},
1 - \delta  \leq y \leq t  \right\}
\end{equation*}
where $C$ is a large positive constant. 
After a change of variables the contribution $J_{\Omega}$ from 
$\sqrt{2n} \cdot \Omega$ to the variance can be written as
\begin{equation*}
J_{\Omega} = \iint_{\Omega} 
\left[ \sqrt{2n} K_n(\sqrt{2n}x, \sqrt{2n}y ) \right]^2 \, \ud x \, \ud y.
\end{equation*}
In order to deal with this integral we must first study the integrand
and, via Christoffel-Darboux, especially the difference
\begin{align} \label{diff}
D = \Big[ h_n \left( \sqrt{2n} x \right) & h_{n-1} \left(\sqrt{2n} y \right) \\
& \qquad \qquad 
- h_{n-1} \left( \sqrt{2n} x \right) h_n \left(\sqrt{2n} y \right) \Big]
e^{-n (x^2 + y^2)}. \nonumber
\end{align}
We will show that in $\Omega$ it holds that
\begin{align*}
b (4n(n-1))^{1/4} D  =  
\bigg[ & \ai \left( - \left[ 3nF(x) \right]^{2/3} \right)
\dai \left( - \left[ 3nF(y) \right]^{2/3} \right) \\ - &  
\dai \left( - \left[ 3nF(x) \right]^{2/3} \right)
\ai \left( - \left[ 3nF(y) \right]^{2/3} \right) \bigg]  \\ 
& \, \, \, \, \, \, \, \, \, \, \, \, \, \, \, \, \, \, \, \, \, \, \,
+ \mathcal{O} \left( \frac{1}{n (1-x)} \right) + 
\mathcal{O} \left( \frac{(1-y)^{3/4}}{(1-x)^{1/4}} \right) .
\end{align*}
Here $b$ is a constant, $\ai$ stands for the Airy function and
\begin{equation*}
F(x) = \int_x^1{\sqrt{1 - t^2}} \, \ud t .
\end{equation*}
In $\Omega$ $h_n$ has the following asymptotics:
\begin{align*}
& h_n(\sqrt{2n} x) e^{-n x^2} \\ & \quad  =
(2n)^{-1/4} \bigg\{ 
\left( \frac{1+x}{1-x} \right)^{1/4} [3nF(x)]^{1/6} 
\ai \left( - [3nF(x)]^{2/3} \right) ( 1 + \mathcal{O} (n^{-1}) ) \\
& \qquad \qquad \quad \, \, \, \, 
- \left( \frac{1-x}{1+x} \right)^{1/4} [3nF(x)]^{-1/6} 
\dai \left( - [3nF(x)]^{2/3} \right) ( 1 + \mathcal{O} (n^{-1}) ) 
\bigg\} 
\end{align*}
If we disregard the $\mathcal{O} (n^{-1})$ terms in
the $h_n-$asymtotics, then (\ref{diff}) can be written as a sum of four 
differences $D_1 - D_4$:
\begin{align*}
(4n(n-1))^{1/4} & D_1 \\
& = \left( \frac{1+x}{1-x} \right)^{1/4} 
\left( \frac{1+y_n}{1-y_n} \right)^{1/4} 
\left[ 3nF(x) \right]^{1/6}  \left[ 3n'F(y_n) \right]^{1/6} \\  
& \,\,\, \,\,\, \,\,\, \,\,\, \,\,\, \,\,\, \,\,\, \,\,\, \,\,\,\,\,\, \,\,\, 
\,\,\,  \,\,\, 
\times \ai \left( - \left[ 3nF(x) \right]^{2/3} \right) 
\ai \left( - \left[ 3n'F(y_n) \right]^{2/3} \right)  \\ & -
\left( \frac{1+x_n}{1-x_n} \right)^{1/4} 
\left( \frac{1+y}{1-y} \right)^{1/4} 
\left[ 3n'F(x_n) \right]^{1/6}  \left[ 3nF(y) \right]^{1/6}  \\
& \,\,\, \,\,\, \,\,\, \,\,\, \,\,\, \,\,\, \,\,\, \,\,\, \,\,\,\,\,\, \,\,\, 
\,\,\,  \,\,\, 
\times \ai \left( - \left[ 3n'F(x_n) \right]^{2/3} \right) 
\ai \left( - \left[ 3nF(y) \right]^{2/3} \right)
\end{align*}
\begin{align*}
(4n(n-1))^{1/4} &  D_2 \\ & = 
\left( \frac{1+x}{1-x} \right)^{1/4} 
\left( \frac{1-y_n}{1+y_n} \right)^{1/4} 
\left[ 3nF(x) \right]^{1/6}  \left[ 3n'F(y_n) \right]^{-1/6} \\
 & \,\,\, \,\,\, \,\,\, \,\,\, \,\,\, \,\,\, \,\,\, \,\,\, \,\,\,\,\,\, \,\,\, 
\,\,\,  \,\,\, 
\times \ai \left( - \left[ 3nF(x) \right]^{2/3} \right) 
\dai \left( - \left[ 3n'F(y_n) \right]^{2/3} \right) \\ & -
\left( \frac{1+x_n}{1-x_n} \right)^{1/4} 
\left( \frac{1-y}{1+y} \right)^{1/4} 
\left[ 3n'F(x_n) \right]^{1/6}  \left[ 3nF(y) \right]^{-1/6} \\
& \,\,\, \,\,\, \,\,\, \,\,\, \,\,\, \,\,\, \,\,\, \,\,\, \,\,\,\,\,\, \,\,\, 
\,\,\,  \,\,\, 
\times \ai \left( - \left[ 3n'F(x_n) \right]^{2/3} \right) 
\dai \left( - \left[ 3nF(y) \right]^{2/3} \right)
\end{align*}
\begin{align*}
(4n(n-1))^{1/4} & D_3 \\ & = 
\left( \frac{1-x_n}{1+x_n} \right)^{1/4} 
\left( \frac{1+y}{1-y} \right)^{1/4} 
\left[ 3n'F(x_n) \right]^{-1/6}  \left[ 3nF(y) \right]^{1/6} \\
& \,\,\, \,\,\, \,\,\, \,\,\, \,\,\, \,\,\, \,\,\, \,\,\, \,\,\,\,\,\, \,\,\, 
\,\,\,  \,\,\,  
\times \dai \left( - \left[ 3n'F(x_n) \right]^{2/3} \right) 
\ai \left( - \left[ 3nF(y) \right]^{2/3} \right) +
\end{align*}
\begin{align*}
& \qquad \qquad \qquad - \left( \frac{1-x}{1+x} \right)^{1/4} 
\left( \frac{1+y_n}{1-y_n} \right)^{1/4} 
\left[ 3nF(x) \right]^{-1/6}  \left[ 3n'F(y_n) \right]^{1/6}  \\
& \qquad \qquad \qquad
\,\,\, \,\,\, \,\,\, \,\,\, \,\,\, \,\,\, \,\,\, \,\,\, \,\,\,\,\,\, \,\,\, 
\,\,\,  \,\,\, 
\times \dai \left( - \left[ 3nF(x) \right]^{2/3} \right) 
\ai \left( - \left[ 3n'F(y_n) \right]^{2/3} \right)
\end{align*}
\begin{align*}
(4n(n-1))^{1/4}&  D_4 \\ & = 
\left( \frac{1-x}{1+x} \right)^{1/4} 
\left( \frac{1-y_n}{1+y_n} \right)^{1/4} 
\left[ 3nF(x) \right]^{-1/6}  \left[ 3n'F(y_n) \right]^{-1/6} \\
& \,\,\, \,\,\, \,\,\, \,\,\, \,\,\, \,\,\, \,\,\, \,\,\, \,\,\,\,\,\, \,\,\, 
\,\,\,  \,\,\,  
\times \dai \left( - \left[ 3nF(x) \right]^{2/3} \right) 
\dai \left( - \left[ 3n'F(y_n) \right]^{2/3} \right) \\  & -
\left( \frac{1-x_n}{1+x_n} \right)^{1/4} 
\left( \frac{1-y}{1+y} \right)^{1/4} 
\left[ 3n'F(x_n) \right]^{-1/6}  \left[ 3nF(y) \right]^{-1/6}  \\
& \,\,\, \,\,\, \,\,\, \,\,\, \,\,\, \,\,\, \,\,\, \,\,\, \,\,\,\,\,\, \,\,\, 
\,\,\,  \,\,\, 
\times \dai \left( - \left[ 3n'F(x_n) \right]^{2/3} \right) 
\dai \left( - \left[ 3nF(y) \right]^{2/3} \right)
\end{align*} 
We have used the notation $n' = n - 1$ and $x_n = \sqrt{\frac{n}{n-1}} x$. 
Note that $x_n < 1$ in $\Omega$.

$D_1$: A calculation using the series expansion 
\begin{equation*}
\frac{F^{1/6}(x)}{(1-x)^{1/4}} =
c_0 + c_1 (1-x) + \ldots
\end{equation*}
gives
\begin{align*}
\left( \frac{1+x_n}{1-x_n} \right)^{1/4} &
\left( \frac{1+y}{1-y} \right)^{1/4} 
\left[ 3(n-1)F(x_n) \right]^{1/6}  \left[ 3nF(y) \right]^{1/6} \\ & =
\left( \frac{1+x}{1-x} \right)^{1/4} 
\left( \frac{1+y}{1-y} \right)^{1/4} 
\left[ 3nF(x) \right]^{1/6}  \left[ 3nF(y) \right]^{1/6} \\ & 
\qquad \qquad \qquad \qquad \qquad  \qquad \qquad \qquad
+ \mathcal{O} \left( n^{1/3} (1-x) \right) \\
& = a_1 n^{1/3} + \mathcal{O} \left( n^{1/3} (1-y) \right) 
\end{align*}
where
\begin{equation*}
a_1 = \lim_{x \rightarrow 1^-} \sqrt{1+x} 
\frac{(3F(x))^{1/3}}{\sqrt{1-x}} .
\end{equation*}
Since 
\begin{equation*}
\ai \left( - \left [3nF(x) \right]^{2/3} \right) =
\mathcal{O} \left( \frac{1}{n^{1/6} (1 - x)^{1/4}} \right)
\end{equation*}
it holds that
\begin{align*}
(4n(n-1))^{1/4} D_1 & \\ = a_1 n^{1/3} 
\bigg[ & \ai \left( - \left[ 3nF(x) \right]^{2/3} \right)
\ai \left( - \left[ 3n'F(y_n) \right]^{2/3} \right)  \\ -  &
\ai \left( - \left[ 3n'F(x_n) \right]^{2/3} \right)
\ai \left( - \left[ 3nF(y) \right]^{2/3} \right) \bigg] \\
& \, \, \, \, \, \, \, \, \, \, \, \, \, \, \, \, \, \, \, \, \, 
\, \, \, \, \, \, \, \, \, \, \, \, \, \, \, \, \, \, \, \, \, 
\, \, \, \, \, \, \, \, \, \, \, \, \, \, \, \, \, \, \, \, \, 
\, \, \, \, \, \, \, \, \, \, \, \, \, \, 
+ \mathcal{O} \left( \frac{(1-y)^{3/4}}{(1-x)^{1/4}} \right) .
\end{align*}
$D_2 - D_4$: The same procedure as above gives
\begin{align*}
(4n(n-1))^{1/4} D_2 & \\ = \mathcal{O} (1) 
\bigg[ & \ai \left( - \left[ 3nF(x) \right]^{2/3} \right)
\dai \left( - \left[ 3n'F(y_n) \right]^{2/3} \right) \\ - &  
\ai \left( - \left[ 3n'F(x_n) \right]^{2/3} \right)
\dai \left( - \left[ 3nF(y) \right]^{2/3} \right) \bigg] \\
&\, \, \, \, \, \, \, \, \, \, \, \, \, \, \, \, \, \, \, \, \, 
\, \, \, \, \, \, \, \, \, \, \, \, \, \, \, \, \, \, \, \, \, 
\, \, \, \, \, \, \, \, \, \, \, \, \, \, \, \, \, \, \, \, \, 
\, \, \, \, \, \, \, \, \, \, \, \, \, \,
+ \mathcal{O} \left( \frac{(1-y)^{5/4}}{(1-x)^{1/4} } \right) 
\end{align*} 
\begin{align*}
(4n(n-1))^{1/4} D_3 & \\  = \mathcal{O} (1) 
\bigg[ & \dai \left( - \left[ 3n'F(x_n) \right]^{2/3} \right)
\ai \left( - \left[ 3nF(y) \right]^{2/3} \right) \\ - &  
\dai \left( - \left[ 3nF(x) \right]^{2/3} \right)
\ai \left( - \left[ 3n'F(y_n) \right]^{2/3} \right) \bigg] \\
&\, \, \, \, \, \, \, \, \, \, \, \, \, \, \, \, \, \, \, \, \, 
\, \, \, \, \, \, \, \, \, \, \, \, \, \, \, \, \, \, \, \, \, 
\, \, \, \, \, \, \, \, \, \, \, \, \, \, \, \, \, \, \, \, \, 
\, \, \, \, \, \, \, \, \, \, \, \, \, \,
+ \mathcal{O} \left( \frac{(1-y)^{5/4}}{(1-x)^{1/4} } \right) 
\end{align*}
\begin{align*}
(4n(n-1))^{1/4} D_4 & \\ = \mathcal{O} (n^{-1/3}) 
\bigg[ & \dai \left( - \left[ 3nF(x) \right]^{2/3} \right)
\dai \left( - \left[ 3n'F(y_n) \right]^{2/3} \right) \\ - &  
\dai \left( - \left[ 3n'F(x_n) \right]^{2/3} \right)
\dai \left( - \left[ 3nF(y) \right]^{2/3} \right) \bigg] \\
&\, \, \, \, \, \, \, \, \, \, \, \, \, \, \, \, \, \, \, \, \, 
\, \, \, \, \, \, \, \, \, \, \, \, \, \, \, \, \, \, \, \, \, 
\, \, \, \, \, \, \, \, \, \, \, \, \, \, \, \, \, \, \, \, \, 
\, \, \, \, \, \, \, \, \, \, \, \, \, \,
+ \mathcal{O} \left( (1-y)^{3/2} \right) 
\end{align*}
Now consider the difference still left in $D_1$:
\begin{multline*}
\ai \left( - \left[ 3nF(x) \right]^{2/3} \right)
\ai \left( - \left[ 3n'F(y_n) \right]^{2/3} \right) \\ -   
\ai \left( - \left[ 3n'F(x_n) \right]^{2/3} \right)
\ai \left( - \left[ 3nF(y) \right]^{2/3} \right)  
\end{multline*}
To deal with this expression we shall first investigate the 
argument
\begin{equation*}
\left[ 3n'F(x_n) \right]^{2/3} = 
\left[ 3(n-1)F\left(\sqrt{\frac{n}{n-1}} x\right) \right]^{2/3}.
\end{equation*}
A simple integration shows that
\begin{equation*}
F(x) = \int_x^1{\sqrt{1-t^2} \, \ud t} = 
\frac{1}{2} \left( \arccos x - x \sqrt{1 - x^2} \right).
\end{equation*}
Since 
\begin{equation*}
x_n = \sqrt{\frac{n}{n-1}} x = x + \frac{x}{2(n-1)} +
\mathcal{O} (n^{-2})
\end{equation*}
we get 
\begin{multline*}
F (x_n) =
F(x) + F'(x) \left( \frac{x}{2(n-1)} + \mathcal{O} (n^{-2}) \right)
+ \mathcal{O} ( F''(x) n^{-2}) \\ =
F(x) - \frac{x \sqrt{1 - x^2}}{2 (n-1)} + 
\mathcal{O} \left( \frac{1}{n^2 \sqrt{1-x}} \right)
\end{multline*}
and hence
\begin{multline*}
3n'F(x_n) = 3(n-1) F(x) - \frac{3}{2} x \sqrt{1-x^2} +
\mathcal{O} \left( \frac{1}{n \sqrt{1-x} } \right) \\ =
3n F(x) - \frac{3}{2} \arccos x +
\mathcal{O} \left( \frac{1}{n \sqrt{1-x} } \right) .
\end{multline*} 
The argument can now finally be written as
\begin{equation} \label{arg}
- \left[ 3n'F(x_n) \right]^{2/3} =
- \left[ 3nF(x) \right]^{2/3} + \frac{\arccos x}{(3nF(x))^{1/3}} +
\mathcal{O} \left( \frac{1}{n^{4/3} (1-x)} \right).
\end{equation}
Note in the last expression that
\begin{equation*}
\frac{\arccos x}{(3nF(x))^{1/3}} \sim n^{-1/3}.
\end{equation*}
It is now possible to expand the difference in a Taylor series around
the point $- \left[ 3nF(x) \right]^{2/3}$ and the result is
\begin{multline*}
\frac{a_2}{n^{1/3}} 
\Big[ \ai \left( - \left[ 3nF(x) \right]^{2/3} \right)
\dai \left( - \left[ 3nF(y) \right]^{2/3} \right) \\ -
\dai \left( - \left[ 3nF(x) \right]^{2/3} \right)
\ai \left( - \left[ 3nF(y) \right]^{2/3} \right) \Big] \\ +
\mathcal{O} \left( \frac{1}{n^{4/3} (1-x)} \right) +
\mathcal{O} \left( \frac{(1-y)^{3/4}}{n^{1/3}(1-x)^{1/4}} \right) 
\end{multline*}
where $a_2$ is defined by
\begin{equation*}
a_2 = \lim_{x \rightarrow 1^-} \frac{\arccos x}{(3F(x))^{1/3}} .
\end{equation*}
Similar computations can be done in $D_2$ - $D_4$ and one then ends 
up with 
\begin{equation*}
(4n(n-1))^{1/4} (D_2 + D_3 + D_4) = 
\mathcal{O} \left( (1-y)^{1/2} \right) .
\end{equation*} 
Adding everything up we now have that
\begin{align} \nonumber 
\frac{(4n(n-1))^{1/4}}{a_1 a_2} D =   
\bigg[ & \ai \left( - \left[ 3nF(x) \right]^{2/3} \right)
\dai \left( - \left[ 3nF(y) \right]^{2/3} \right) \\ - & \label{D} 
\dai \left( - \left[ 3nF(x) \right]^{2/3} \right)
\ai \left( - \left[ 3nF(y) \right]^{2/3} \right) \bigg]  \\ \nonumber
& \, \, \, \, \, \, \, \, \, \, \, \, \, \, \, \, \, \, \, \, \, \, \,
\, \, \, \, \, \, \, +
\mathcal{O} \left( \frac{1}{n (1-x)} \right) + 
\mathcal{O} \left( \frac{(1-y)^{3/4}}{(1-x)^{1/4}} \right) .
\end{align} 
As we shall see the main contribution comes from the domain 
\begin{equation*}
\Omega_1 = \left\{ (x,y) ; t \leq x \leq t + \frac{1-t}{r(n)} , 
t - \frac{1-t}{r(n)} \leq y \leq t - \epsilon \right\} .
\end{equation*}
Here $r(n)$ is a function tending slowly to infinity as $n$ tends to 
infinity and $\epsilon = \frac{1}{n(1-t)^{1/2}}$, the size of the expected 
distance between two eigenvalues at $t$. The reason why this 
$\epsilon$ is necessary lies in the asymptotics for the Hermite
polynomials. The error term given there, however small, will cause 
problems since the integral
\begin{equation*}
\int_{t}^{t + \epsilon} {
\int_{t - \epsilon}^{t}{ \frac{1}{(x-y)^2} \, \udy } \, \udx }
\end{equation*}
is divergent.

From the asymptotics of the Airy function and its 
derivative we get, for $(x,y) \in \Omega_1$, that  
\begin{align*}
\ai & \left( - \left[ 3nF(x) \right]^{2/3} \right)
\dai \left( - \left[ 3nF(y) \right]^{2/3} \right) \\ & =
\left( \frac{1}{\sqrt{\pi}} (nF(x))^{-1/6}
\sin{ \left[ 2nF(x) + \frac{\pi}{4} \right]} +
\mathcal{O} \left( (nF(x))^{-7/6} \right) \right) \\ & \qquad \qquad \times
\left( \frac{1}{\sqrt{\pi}} (nF(y))^{1/6}
\sin{ \left[ 2nF(y) - \frac{\pi}{4} \right] } +
\mathcal{O} \left( (nF(y))^{-5/6} \right) \right) \\ & =
\frac{1}{\pi} \left( \frac{F(y)}{F(x)} \right)^{1/6}
\sin{ \left[ 2nF(x) + \frac{\pi}{4} \right] } 
\sin{ \left[ 2nF(y) - \frac{\pi}{4} \right] } + 
\mathcal{O} \left( (nF(x))^{-1} \right) .
\end{align*}
If we define $r(n)$ by
\begin{equation*} 
\frac{1}{r(n)} = \max \left( \sqrt{1-t}, 
\frac{1}{\log [n (1-t)^{3/2}]} \right)
\end{equation*} 
we get that
\begin{align*}
\left( \frac{F(y)}{F(x)} \right)^{1/6} & =
1 + \mathcal{O} \left( \frac{1}{r(n)} \right) \\
\left( \frac{F(x)}{F(y)} \right)^{1/6} & =
1 + \mathcal{O} \left( \frac{1}{r(n)} \right) \\
\frac{1}{n (1-x)} \, \, \, \, \, \, & =
\mathcal{O} \left( \frac{1}{n F(x)} \right) =
\mathcal{O} \left( \frac{1}{r(n)} \right) \\
\frac{(1-y)^{3/4}}{(1-x)^{1/4}}  \, \, \, & = 
\mathcal{O} \left( \frac{1}{r(n)} \right) .
\end{align*}
From this it follows that in $\Omega_1$ $D$ can be written as
\begin{align*}
\frac{(4n(n-1))^{1/4}}{a_1 a_2} & D \\ = &
\frac{1}{\pi} 
\bigg( \sin{ \left[ 2nF(x) + \frac{\pi}{4} \right] }  
\sin{ \left[ 2nF(y) - \frac{\pi}{4} \right] } \\ & \, \, \,   -
\sin{ \left[ 2nF(x) - \frac{\pi}{4} \right] }  
\sin{ \left[ 2nF(y) + \frac{\pi}{4} \right] } \bigg) 
+ \mathcal{O} \left( \frac{1}{r(n)} \right) \\
= & \frac{1}{\pi} 
\sin{ \left[ 2n (F(x) - F(y)) \right] } + 
\mathcal{O}\left( \frac{1}{r(n)} \right).
\end{align*}
The numerator in the integral of interest is 
\begin{multline*}
\frac{n}{2\sqrt{n(n-1)}} D^2 \\ =
\frac{(a_1 a_2)^2}{4 \pi^2} 
\sin^2{ \left[ 2n(F(x)-F(y)) \right]} + 
\mathcal{O}\left( \frac{1}{r(n)} \right) \\ =
\frac{1}{2 \pi^2} 
(1 - \cos{ \left[ 4n(F(x)-F(y)) \right]}) +
\mathcal{O}\left( \frac{1}{r(n)} \right).
\end{multline*} 
We used here that $a_1 a_2 = 2$. A simple integration gives
\begin{equation*}
\iint_{\Omega_1} \frac{1}{(x-y)^2}\, \ud x \, \ud y =
\log[n(1-t)^{3/2}] + \mathcal{O} (\log{r(n)}) .
\end{equation*}
The integral 
\begin{equation*}
I = \iint_{\Omega_1} 
\frac{\cos{ \left[ 4n(F(x)-F(y)) \right]}}{(x-y)^2}\, \ud x \, \ud y
\end{equation*}
is $\mathcal{O} (1)$:
By doing a partial integration $I$ can be split into two 
integrals: 
\begin{multline*}
I = \int_{t}^{t+\frac{1-t}{r(n)}} \Bigg(
\left[ \frac{\sin \left[ 4n(F(x)-F(y)) \right]}{-4nF'(y) (x-y)^2} 
\right]_{t-\frac{1-t}{r(n)}}^{t-\epsilon} \\ 
 + \int_{t-\frac{1-t}{r(n)}}^{t-\epsilon} 
\sin \left[ 4n(F(x)-F(y)) \right]
\left( \frac{1}{4n F'(y)(x-y)^2} \right)_y' 
\, \udy \Bigg) \, \udx \\
=: I_1 + I_2
\end{multline*}
We can estimate $I_1$ in the following way.
\begin{multline*}
|I_1| \leq 2 \int_{t}^{t+\frac{1-t}{r(n)}}
\frac{1}{4n \sqrt{1-t} (x-(t-\epsilon))^2} \, \udx \\ =
\frac{\epsilon}{2} 
\left[ \frac{-1}{x - t + \epsilon} \right]_t^{t+\frac{1-t}{r(n)}} \leq
\frac{\epsilon}{2} \cdot \frac{2}{\epsilon} = 1
\end{multline*}
Since 
\begin{equation*}
\left( \left[ F'(y)(x-y)^2 \right]^{-1} \right)_y' =
- \frac{y}{(1-y^2)^{3/2} (x-y)^2} - \frac{2}{\sqrt{1-y^2} (x-y)^3}
\end{equation*}
we get
\begin{equation*}
|I_2| \leq C \left( \iint_{\Omega_1} \frac{1}{n(1-y)^{3/2}(x-y)^2} 
\, \udx \udy +
\iint_{\Omega_1} \frac{1}{n \sqrt{1-y} (x-y)^3} \, \udx \udy\right) .
\end{equation*}
The first part is small:
\begin{multline*}
\iint_{\Omega_1} \frac{1}{n(1-y)^{3/2}(x-y)^2} \, \udx \udy 
\\ \leq 
\frac{1}{n (1-t)^{3/2}} 
\iint_{\Omega_1} \frac{1}{(x-y)^2} \, \udx \udy = 
\mathcal{O} \left(
\frac{\log \left[ n (1-t)^{3/2} \right]}{n (1-t)^{3/2}} \right)
\end{multline*}
The second part is also easily estimated:
\begin{equation*}
\iint_{\Omega_1} \frac{1}{n \sqrt{1-y}(x-y)^3} \, \udx \udy 
\leq 
\epsilon \iint_{\Omega_1} \frac{1}{(x-y)^3} \, \udx \udy
= \mathcal{O} (1)
\end{equation*}
This concludes the calculations in $\Omega_1$.

The calculations made above can also be applied to the thin slice
\begin{equation*}
\left\{ (x,y) ; t + \epsilon \leq x \leq t + \frac{1-t}{r(n)} ,
t - \epsilon \leq y \leq t \right\}
\end{equation*}
and the result is $\mathcal{O} (\log [r(n)])$. 

The corner
\begin{equation*}
\Omega_0 = \left\{ (x,y) ; t \leq x \leq t + \epsilon ,
t - \epsilon \leq y \leq t \right\}
\end{equation*}
requires a special technique.
In this domain a different representation of $K_n$ will be
used, namely
\begin{equation*}
K_n (x,y) = \sum_{i=0}^{n-1} p_i (x) p_i (y) 
e^{-\frac{1}{2} \left( x^2 + y^2 \right)}.
\end{equation*} 
By use of the Cauchy-Schwartz inequality we get that
\begin{equation*}
K_n^2 (x,y) \leq K_n (x,x) K_n (y,y) .
\end{equation*}
Having separated the variables we can now use the calculations 
from the proof of Lemma \ref{exp2} to see that
\begin{equation*}
\int_{t-\epsilon}^{t} \int_{t}^{t+\epsilon} 
\left( \sqrt{2n} K_n (\sqrt{2n}x,\sqrt{2n}y) \right)^2 \, \udx \udy
 = \mathcal{O} (1) .
\end{equation*}
Note that
\begin{equation*}
\int_{t}^{t+\epsilon} K_n(\sqrt{2n}x,\sqrt{2n}x) \, \udx = 
g(t) - g(t+\epsilon)
\end{equation*}
where $g(t)$ is the expected number of eigenvalues in the interval 
$(t \sqrt{2n}, \infty)$.

Now we shall look at the other part still left of 
$\Omega$. This domain can conveniently be written as 
$\Omega_2 \cup \Omega_3$ where
\begin{equation*}
\Omega_2 = \left\{ (x,y) ; t \leq x \leq 1 - Cn^{-1} ,
1 - \delta \leq y \leq t - \frac{1-t}{r(n)} \right\}
\end{equation*}
and 
\begin{equation*}
\Omega_3 = \left\{ (x,y) ; t + \frac{1-t}{r(n)} 
\leq x \leq 1 - Cn^{-1}  ,
t - \frac{1-t}{r(n)} \leq y \leq t \right\}
\end{equation*}
When looking at the expression for $D$ in (\ref{D}) above it is 
clear that every term is smaller than 
\begin{equation*}
n^{-1/2} \ai \left( - \left[ 3nF(x) \right]^{2/3} \right) 
\dai \left( - \left[ 3nF(y) \right]^{2/3} \right) =
\mathcal{O} \left( n^{-1/2} \left( \frac{1-y}{1-x} \right)^{1/4} \right) .
\end{equation*}
This means that it is sufficient to calculate the integrals
\begin{equation*}
\iint_{\Omega_i} \frac{\sqrt{1-y}}{\sqrt{1-x} (x-y)^2} \, \udx \, \udy
\quad \mathrm{i = 2,3}.
\end{equation*}
The calculations are straightforward so some details will be skipped.
When first integrating with respect to the x-variable one gets
\begin{multline*}
\int_{L_1}^{H_1} \frac{\sqrt{1-y}}{\sqrt{1-x} (x-y)^2} \, \udx \\ =
\frac{1}{2(1-y)} \log \left[ \frac{\left( \sqrt{1-y}+\sqrt{1-L_1} 
\right) \left( \sqrt{1-y}-\sqrt{1-H_1} \right) }
{\left( \sqrt{1-y}+\sqrt{1-H_1} \right) 
\left( \sqrt{1-y}-\sqrt{1-L_1} \right) } \right] \\ +
\frac{1}{\sqrt{1-y}} \left(
\frac{1}{\sqrt{1-y}-\sqrt{1-L_1}} - \frac{1}{\sqrt{1-y}+\sqrt{1-L_1}}
\right) .
\end{multline*}
$\Omega_2$: Letting $H_1 = 1$ instead of $1-Cn^{-1}$ we get nicer 
expressions. This is allowed since the domain of integration becomes
larger. The task is to get an upper bound for the integrals
\begin{multline*}
A = \int_{L_2}^{H_2} \frac{1}{2(1-y)} 
\log \left[ \frac{\sqrt{1-y} + \sqrt{1-L_1}}{\sqrt{1-y} - \sqrt{1-L_1}} 
\right] \, \udy \\ = 
\int_{\sqrt{1-H_2}}^{\sqrt{1-L_2}} \frac{1}{z} 
\log \left[ \frac{z + \sqrt{1-L_1}}{z - \sqrt{1-L_1}} \right] \, \udz
\end{multline*}
and
\begin{multline*}
B = \int_{L_2}^{H_2} \frac{1}{\sqrt{1-y}} \left(
\frac{1}{\sqrt{1-y}-\sqrt{1-L_1}} - 
\frac{1}{\sqrt{1-y}+\sqrt{1-L_1}} \right) \, \udy \\ =
2 \int_{\sqrt{1-H_2}}^{\sqrt{1-L_2}} \left(
\frac{1}{z-\sqrt{1-L_1}} - \frac{1}{z+\sqrt{1-L_1}} \right) \, \udz
\end{multline*}
where 
\begin{equation*}
L_2 = 1- \delta \, \mathrm{,} \quad H_2 = t - \frac{1-t}{r(n)}
\quad \mathrm{and} \quad L_1 = t \, .
\end{equation*}
When manipulating the integrand in A one gets
\begin{equation*}
\frac{1}{z} \log \left[ 1 + 2 \frac{\sqrt{1-L_1}}{z - \sqrt{1-L_1}}
\right] = \frac{1}{z} \mathcal{O} \left(
\frac{\sqrt{1-L_1}}{z - \sqrt{1-L_1}} \right) \, .
\end{equation*}
Some algebra shows that
\begin{equation*}
\frac{\sqrt{1-L_1}}{ z \left( z - \sqrt{1-L_1} \right)} =
\frac{1}{z - \sqrt{1-L_1} } - \frac{1}{z} 
\end{equation*} 
which can easily be integrated:
\begin{equation*}
A \leq C \left[ \log \left[
\frac{z - \sqrt{1-L_1}}{z} \right] 
\right]_{\sqrt{1-H_2}}^{\sqrt{1-L_2}}  =  
\mathcal{O} (\log r(n))
\end{equation*} 
The integral B is even easier and one gets
\begin{equation*}
B = 2 \left[ \log \left[ 
\frac{z-\sqrt{1-L_1}}{z+\sqrt{1-L_1}} \right] 
\right]_{\sqrt{1-H_2}}^{\sqrt{1-L_2}} = \mathcal{O} (\log r(n)) \, .
\end{equation*}
\noindent $\Omega_3$: The same procedure as in $\Omega_2$ gives that the
contribution to the variance from this domain is $o(1)$.

We shall now consider the thin strip
\begin{equation*}
\Omega_4 = \left \{ x,y ; 1 - C n^{-1} \leq x \leq 1 + C n^{-1}, 
1 - \delta \leq y \leq t \right\}.
\end{equation*}
The asymptotics here is similar to that in $\Omega$ and hence many of the
calculations already done can be applied here as well. As before $D$ can
be split up in $D_1 - D_4$ which can all be treated similarly. Therefore
we only look at $D_1$ here. We have that
\begin{align*}
(4n(n-1))^{1/4} D_1 = a_1 n^{1/3} 
\bigg[ & \ai \left( \mp \left[ 3nF(x) \right]^{2/3} \right)
\ai \left( - \left[ 3n'F(y_n) \right]^{2/3} \right) \\ - & 
\ai \left( \mp \left[ 3n'F(x_n) \right]^{2/3} \right)
\ai \left( - \left[ 3nF(y) \right]^{2/3} \right) \bigg] \\ 
& \,\,\,\,\, \,\,\,\,\, \,\,\,\,\, \,\,\,\,\, \,\,\,\,\, \,\,\,\,\, 
\,\,\,\,\, \,\,\,\,\, \,\,\,\,\, \,\,\,\,\, 
\,\,\,\,\, \,\,\,\,\, \,\,\,\,\, \,\,\,\,\, 
+ \mathcal{O} \left( \frac{(1-y)^{3/4}}{(1-x)^{1/4}} \right). 
\end{align*}
where 
\begin{equation*}
\mp \left[ 3nF(x) \right]^{2/3} = \bigg\{
\begin{array}{cl}
- \left[ 3nF(x) \right]^{2/3} & \quad \textrm{if } x < 1 \\
\left[ 3nF(x) \right]^{2/3} & \quad \textrm{if } x \geq 1 .
\end{array}
\end{equation*}
This follows from the calculations done
above and the asymptotics for the Hermite polynomials when $x>1$. 
In $\Omega_4$ we have
\begin{eqnarray*}
\ai \left( \mp \left[ 3nF(x) \right]^{2/3} \right) &=& 
\ai(0) + \mathcal{O}(n^{-1/3}) \\
\ai \left( \mp \left[ 3n'F(x_n) \right]^{2/3} \right) &=& 
\ai(0) + \mathcal{O}(n^{-1/3})
\end{eqnarray*}
and by using equation (\ref{arg}) (for the $y$-variable) one gets
\begin{equation*}
(4n(n-1))^{1/4} D_1 =
\mathcal{O} \left( \frac{(1-y)^{1/4}}{|1-x|^{1/4}} \right).
\end{equation*}
The error term here has actually already been dealt with in the 
estimations of the contribution coming from $\Omega_2$. 

Rather than to repeat a lot of calculations we now just give ideas
of how to to treat what is left of $[t,\infty) \times (-\infty,t]$.

In the domain
\begin{equation*}
\left\{ x,y ; 1 + C n^{-1} \leq x \leq 1 + \delta, 1- \delta \leq y \leq t 
\right\}
\end{equation*}
one can perform much the same calculations as in $\Omega$ and the 
contribution is $\mathcal{O}(1)$. In
\begin{equation*}
\{ x,y ; t \leq x \leq 1 + \delta, -1 - \delta \leq y \leq 1 - \delta \}
\end{equation*}
one can use the fact that $x-y \geq \delta$ to show
that the contribution from this domain is $\mathcal{O}(1)$.
If $x \geq 1 + \delta$ or $y \leq -1 - \delta$ t
one easily gets from the asymptotics for the Hermite
polynomials that $K_n(\sqrt{2n}x,\sqrt{2n}y)$ is exponentially small 
in $n$ and exponentially
decaying in $x^2$ (or $y^2$). Thus the contribution from this domain
is $o(1)$.
\end{proof}

\section*{appendix}
\setcounter{equation}{0}
\subsection*{A: Some integrals.}

The following equalities hold:
\begin{eqnarray*}
\int_x^{\infty} \ai^2(y) \udy &=&
\ai^{{\prime}^2}(x) - x \ai^2(x) \\
\int_x^{\infty} y \ai^2(y) \udy &=&
\frac{1}{3} 
\left( x \ai^{{\prime}^2}(x) - x^2 \ai^2(x) - \ai(x)\dai(x) \right) \\
\int_x^{\infty} \ai^{{\prime}^2}(y) \udy &=& 
\frac{1}{3} 
\left( x^2 \ai^2(x) - x \ai^{{\prime}^2}(x) - 2 \ai(x)\dai(x) \right) \\
\int_x^{\infty} y^2 \ai^2(y) \udy &=&
\frac{1}{5} \left( x^2 \ai^{{\prime}^2}(x) - x^3 \ai^2(x) - 2x \ai(x)\dai(x) 
+ \ai^2(x) \right) \\
\int_x^{\infty} y \ai^{{\prime}^2}(y) \udy &=&
\frac{1}{5} \left( x^3 \ai^2(x) - x^2 \ai^{{\prime}^2}(x) - 3x \ai(x)\dai(x) 
+ \frac{3}{2} \ai^2(x) \right)
\end{eqnarray*}
The first integral is obtained by performing a partial integration while 
remembering that 
\begin{equation*}
\textrm{Ai}''(x) = x \ai(x).
\end{equation*}
The integrals 3-5 can be obtained rather easily from the second which 
can be treated as follows: \\
Set
\begin{equation*}
u_{\alpha} (x) = \ai(\alpha x) \quad \alpha > 0 .
\end{equation*}
The relationship 
\begin{equation*}
\left[ u_{\alpha}' u_{\beta} - u_{\alpha} u_{\beta}' \right]' =
u_{\alpha}'' u_{\beta} - u_{\alpha} u_{\beta}'' =
x (\alpha^3 - \beta^3) u_{\alpha} u_{\beta} 
\end{equation*}
holds since 
\begin{equation*}
u_{\alpha}'' (x) = \alpha^2 \textrm{Ai}''(\alpha x) = 
\alpha^3 x \ai(\alpha x) = \alpha^3 x u_{\alpha} (x).
\end{equation*} 
Hence 
\begin{multline*}
\int_a^{\infty} x u_{\alpha}(x) u_{\beta}(x) \, \udx =
\frac{1}{\alpha^3 - \beta^3} 
\left[ u_{\alpha}' u_{\beta} - u_{\alpha} u_{\beta}' 
\right]_a^{\infty} \\ =
\frac{u_{\alpha}(a) u_{\beta}'(a) - 
u_{\alpha}'(a) u_{\beta}(a)}{\alpha^3 - \beta^3} .
\end{multline*}
The idea now is to let $\alpha$, $\beta$ tend to one.
Set $\alpha = 1+h$ and $\beta = 1-h$ where $h>0$ and small. 
The left hand side tends to
\begin{equation*}
\int_a^{\infty} x \ai^2 (x) \, \udx
\end{equation*}
as $h \rightarrow 0^+$. Standard calculations show that 
at the same time the right hand side tends to
\begin{equation*}
\frac{1}{3} \left(
- a^2 \ai^2(a) - \ai(a) \dai(a) + a \ai^{{\prime}^2}(a) \right) .
\end{equation*}

\subsection*{B: Variance calculations.}
Let $I_1, \ldots, I_m$ be a set of disjoint intervals and 
$\#I_i$ be the number of eigenvalues of the $\textrm{GUE}_n$ 
in the interval $I_i$. We shall give a formula for 
$\var \left( \alpha_1 \#I_1 + \ldots + \alpha_m \#I_m\right)$.
Clearly
\begin{equation*}
\#I_i = \sum_{k=1}^n \chi_{I_i}(x_k) \qquad 1 \leq i \leq n
\end{equation*}
where $\chi_B$ is the characteristic function for the set $B$ and 
$\{x_k\}_{k=1}^n$ are the unordered eigenvalues.  
The expected value is easy to compute:
\begin{equation*}
\E \left[ \#I_i \right] =
\int_{I_i} \rho_{n,1} (x) \, \udx = \int_{I_i} K_n(x,x) \, \udx
\end{equation*}
The correlation functions $\rho_{n,k}$ were defined in the introduction. 
We also need to calculate $\E[\#I_i^2]$:
\begin{align*}
\E \left[ \#I_i^2 \right] &=
\E \left[ \sum_{j,k=1}^n  \chi_{I_i}(x_k) \chi_{I_i}(x_j)\right] \\ &=
\sum_{k=1}^n  \E [ \chi_{I_i}(x_k) ] +
\sum_{j \neq k}  \E [\chi_{I_i}(x_k) \chi_{I_i}(x_j) ] \\ & =
\int_{I_i} K_n(x,x) \, \udx +
\iint_{I_i \times I_i} \rho_{n,2} (x,y) \, \udx \udy \\ & =
\int_{I_i} K_n(x,x) \, \udx +
\left( \int_{I_i} K_n(x,x) \, \udx \right)^2 -
\iint_{I_i \times I_i} K_n^2(x,y) \, \udx \udy
\end{align*}
It now follows that
\begin{equation*}
\var(\#I_i) = \int_{I_i} K_n(x,x) \, \udx -
\iint_{I_i \times I_i} K_n^2(x,y) \, \udx \udy .
\end{equation*}
To get a more convenient formula to work with one can now use the
identities $K_n(x,y) = K_n(y,x)$ and
\begin{equation*}
\int_{\mathbb{R}} K_n(x,y) K_n(y,z) \, \udy = K(x,z)
\end{equation*}
to get
\begin{multline*}
\var(\#I_i) = \int_{I_i} \left(  
\int_{\mathbb{R}} K_n^2(x,y) \, \udy \right) \, \udx -
\iint_{I_i \times I_i} K_n^2(x,y) \, \udx \udy \\ =
\iint_{I_i \times I_i^c} K_n^2(x,y) \, \udx \udy .
\end{multline*}
In more generality one gets
\begin{equation*}
\E [ \alpha_1 \#I_1 + \ldots + \alpha_m \#I_m  ] =
\sum_{i=1}^m \alpha_i \int_{I_i} K_n(x,x) \, \udx 
\end{equation*} 
and 
\begin{multline*}
(\alpha_1 \#I_1 + \ldots + \alpha_m \#I_m)^2 \\ =
\sum_{i=1}^m \alpha_i^2 \left( \sum_{k=1}^n \chi_{I_i}(x_k) \right)^2 
+ \sum_{i \neq j}^m \alpha_i \alpha_j
\left( \sum_{k=1}^n \chi_{I_i}(x_k) \right)
\left( \sum_{k=1}^n \chi_{I_j}(x_k) \right) \\ =: S_1 + S_2.
\end{multline*}
From the calculations above we know that
\begin{multline*}
\E[S_1] = \sum_{i=1}^m
\alpha_i^2 \bigg\{  \int_{I_i} K_n(x,x) \, \udx +
\left( \int_{I_i} K_n(x,x) \, \udx \right)^2 \\ -
\iint_{I_i \times I_i} K_n^2(x,y) \, \udx \udy \bigg\}
\end{multline*}
so it remains to calculate $\E[S_2]$. We have 
\begin{equation*}
\left( \sum_{k=1}^n \chi_{I_i}(x_k) \right)
\left( \sum_{k=1}^n \chi_{I_j}(x_k) \right) =
\sum_{k \neq l}^n \chi_{I_i}(x_k) \chi_{I_i}(x_l)
\end {equation*}
and hence 
\begin{multline*}
\E[S_2] = \sum_{i \neq j}^m \alpha_i \alpha_j
\iint_{I_i \times I_i} \rho_{n,2} (x,y) \, \udx \udy \\ =
\sum_{i \neq j}^m \alpha_i \alpha_j \left(
\int_{I_i} K_n(x,x) \, \udx \int_{I_j} K_n(x,x) \, \udx
- \iint_{I_i \times I_j} K_n^2(x,y) \, \udx \udy \right).
\end{multline*}
Since 
\begin{multline*}
\left( \E[\alpha_1 \#I_1 + \ldots + \alpha_m \#I_m ] \right)^2 \\ =
\sum_{i=1}^m \alpha_i^2 \left( \int_{I_i} K_n(x,x) \, \udx \right)^2 + 
\sum_{i \neq j}^m \alpha_i \alpha_j 
\int_{I_i} K_n(x,x) \, \udx \int_{I_j} K_n(x,x) \, \udx
\end{multline*}
we finally get (with manipulations as before)
\begin{multline*}
\var (\alpha_1 \#I_1 + \ldots + \alpha_m \#I_m) \\ =
\sum_{i=1}^m \alpha_i^2 
\iint_{I_i \times I_i^c} K_n^2(x,y) \, \udx \udy -
\sum_{i \neq j}^m \alpha_i \alpha_j 
\iint_{I_i \times I_j} K_n^2(x,y) \, \udx \udy .
\end{multline*}

\section*{Acknowledgment}
I would like to thank my advisor Kurt Johansson 
for presenting this problem to me and for many useful discussions. 
I also thank Kenneth McLaughlin for drawing my attention to his and 
N. M. Ercolani's work on the asymptotics of the Partition Function for 
random matrices.

\end{document}